\documentclass[pdflatex,sn-mathphys-num]{sn-jnl}

\usepackage{graphicx}
\usepackage{multirow}
\usepackage{amsmath,amssymb,amsfonts}
\usepackage{amsthm}
\usepackage{mathrsfs}
\usepackage[title]{appendix}
\usepackage[dvipsnames]{xcolor}
\usepackage{textcomp}
\usepackage{manyfoot}
\usepackage{booktabs}
\usepackage{algorithm}
\usepackage{algorithmicx}
\usepackage{algpseudocode}
\usepackage{listings}

\theoremstyle{thmstyleone}
\newtheorem{theorem}{Theorem}
\newtheorem{proposition}[theorem]{Proposition}
\theoremstyle{thmstyletwo}

\newtheorem{remark}{Remark}
\theoremstyle{thmstylethree}

\usepackage{comment}
\usepackage{tikz}
\usepackage{circuitikz}

\usetikzlibrary{decorations.pathreplacing, calligraphy}

\definecolor{darkblue}{rgb}{0,0,190}   
\definecolor{darkgray}{rgb}{0.55,0.55,0.55}
\definecolor{lightgray}{rgb}{0.827,0.827,0.827}
\usepackage{standalone}

\newcommand{\opinion}[1]{\tikz[baseline]{\draw[fill=#1,line width=0.5pt] (0,0.1) circle(0.9ex)}}

\newcommand{\oO}{\mathcal{O}}

\DeclareMathOperator*{\argmax}{arg\,max}

\begin{document}

\title{Speaking of Opinions: Comparing Approaches to Modelling Opinion Manipulation}

\author*[1]{\fnm{Luisa} \sur{Estrada}}\email{Luisa-Fernanda.Estrada-Plata@warwick.ac.uk}
\equalcont{These authors contributed equally to this work.}

\author[1]{\fnm{Sasha} \sur{Glendinning}}\email{Sasha.Glendinning@warwick.ac.uk}
\equalcont{These authors contributed equally to this work.}

\author[1,2]{\fnm{Andrew} \sur{Nugent}}\email{Andrew.Nugent@ucl.ac.uk}

\affil[1]{\orgdiv{Department of Mathematics}, \orgname{University of Warwick}, 
\country{United Kingdom}}
\affil[2]{\orgdiv{Department of Mathematics}, \orgname{University College London}, 
\country{United Kingdom}}

\abstract{This review outlines the major approaches to modelling opinion formation and manipulation in mathematics and computer science. Key tools such as ordinary and partial differential equations, stochastic models, control theory, and interaction protocols are introduced and compared as methods for describing manipulation. The review is separated into those models using a continuous opinion space and those using discrete or binary opinions, with the advantages and disadvantages of each discussed. Finally, the authors provide an interdisciplinary perspective on the field of opinion dynamics and its social significance.}


\maketitle

\section{Introduction} \label{Section: Introduction}

Models of opinion formation, in which a group of individuals interact and update their beliefs, have been studied over the years from multiple perspectives. Foundational models, such as those of French \cite{french1956formal}, Axelrod \cite{axelrod1997dissemination} and DeGroot \cite{degroot1974reaching}, have given way to a vast literature in both mathematics and computer science. A wide array of modelling techniques have been used, including: game theory \cite{CS404_guideBook, kearns_introduction_1994, etesami_game-theoretic_2015}, agent-based models \cite{hegselmann2015opinion,deffuant2000mixing,nugent2024bridging}, ordinary and partial differential equations \cite{ceragioli2012continuous,motsch2014heterophilious,goddard2022noisy}, and tools from statistical physics \cite{toscani2006kinetic,pareschi2013interacting}. Each approach captures the mechanisms of opinion formation slightly differently and naturally lends itself to addressing different questions.

In this review, we focus on models of opinion \textit{manipulation} from the perspectives of continuous mathematical models of opinion formation and discrete models, commonly studied in computer science. The modelling of opinion manipulation has become an important topic in recent years due to the advent of social media where bots and algorithmic social network manipulation have caused interference in election results, for example \cite{woolley2020bots}. Through this review, we study the similarities and differences between the ways in which the mathematical and computer science communities have addressed modelling manipulation. This includes differences in how well continuous and discrete models capture effects of manipulation and the impact that this has on both research questions and results.

A continuous opinion dynamics model is defined such that individuals in a population have an opinion lying in a continuous space. This allows for rich dynamics on opinion space as individuals update their opinions continuously according to the result of interactions within the population. Continuous opinion dynamics models are generally written as ordinary differential equations or partial differential equations but the field is also strongly linked to agent-based models of opinion formation \cite{nugent2024bridging}. Manipulation is typically modelled by combining the field of continuous opinion dynamics with optimal control methods. In the majority of the models presented here, the individuals in a population are influenced by an external figure seeking to cause a change in opinions \cite{albi2014kinetic,albi2014boltzmann,nugent2024steering} but we also consider a case where the figure driving manipulation is itself an agent in the system \cite{glendinning2025what}.

In contrast, binary opinion dynamic models make the assumption that each individual in a population has one of two opinions, often interpreted as being `yes' or `no'. Many models employ stochastic update rules, such as the classical voter model. In this framework, manipulation is often modelled by introducing interaction protocols, where researchers design an algorithm to complete a given task and then ask if this algorithm succeeds in finite time, what the cost of said algorithm would be and whether we can improve on this \cite{ninjas2018}. Similarly to in the continuous regime, we can also view some binary interaction models of opinion manipulation as optimisation problems, often asking how we can maximise the influence that a stubborn individual, sometimes known as a controller, has over a population \cite{kempe_maximizing_2003,brede2018resisting,kozitsin2022general}.

In comparing approaches to modelling manipulation in mathematics and theoretical computer science, we bring together two perspectives on a similar problem from the mathematical sciences. We explore synergies and contrasts between the two fields and discuss the limitations and strengths of each perspective. As early-career researchers from communities that have historically been underrepresented in the mathematical sciences, we would like to highlight the strength that comes from embracing different perspectives, both mathematically and in terms of fostering an inclusive and diverse community.

Our review is structured as follows. In Section \ref{Section: Cts opinion dynamics}, we introduce models of continuous opinion dynamics before discussing four different frameworks for modelling manipulation. We conclude this section with a discussion. In Section \ref{Section: Binary Opinion Dynamics}, we introduce models of binary opinion formation and discuss game theoretic as well as optimisation problems that arise from considering tactics of manipulation. We also consider here how we can perform inference on the state of a social network by manipulating the opinion dynamics. This section is concluded with a discussion. Finally, Section \ref{Section: conclusion} brings together the two approaches considered in this review and gives our perspectives on the ethics and history of studying manipulation.

\section{Continuous Opinion Dynamics} \label{Section: Cts opinion dynamics}

In this section, we introduce models of opinion formation with a continuous opinion space and describe how opinion manipulation can be approached from a dynamical systems perspective as the introduction of a control to an existing model of opinion dynamics. Specifically, we consider models that are continuous in both time and the opinion space, which will enable the use of optimal control theory. This class of models includes ordinary differential equations (ODEs) describing the evolution of agents’ opinions and partial differential equations (PDEs) describing the evolution of population-level opinion distributions. Models tracking agents' opinions individually are referred to as microscopic models, and those tracking the population-level opinion distribution are referred to as macroscopic models. 

In this context, a natural way to model manipulation is to introduce a control variable into the dynamical system. The first control we consider directly contributes to the evolution equations for individuals' opinions, modelling influence on opinions through advertising, for example. Next we split the population into two groups: leaders and followers and show how the opinion distribution of followers can be shaped by the leaders. Finally we consider a control that does not affect opinions directly but instead affects the social network and thus the way that individuals interact. These controls increase in complexity and we will demonstrate a variety of techniques for showing the existence of and finding controls. In all three cases the goal of controls is to bring the population to consensus as some specified target opinion. 

\subsection{Continuous Opinion Dynamics Models}\label{Section: Introduction continuous models}

Denote by $\mathcal{I}$ an interval in which an agent's opinion, denoted $x_i$, is assumed lie. A common choice is to take $\mathcal{I} = [-1,1]$ where we can interpret an opinion $x_i = 1$ to be an extreme positive opinion and $x_i = -1$ to be an extreme negative opinion on a given topic. Alternatively, we may consider this interval as representing the strength of (dis)agreement with a given statement. 

\begin{figure}[ht!]
    \centering
    \includegraphics[width=\linewidth]{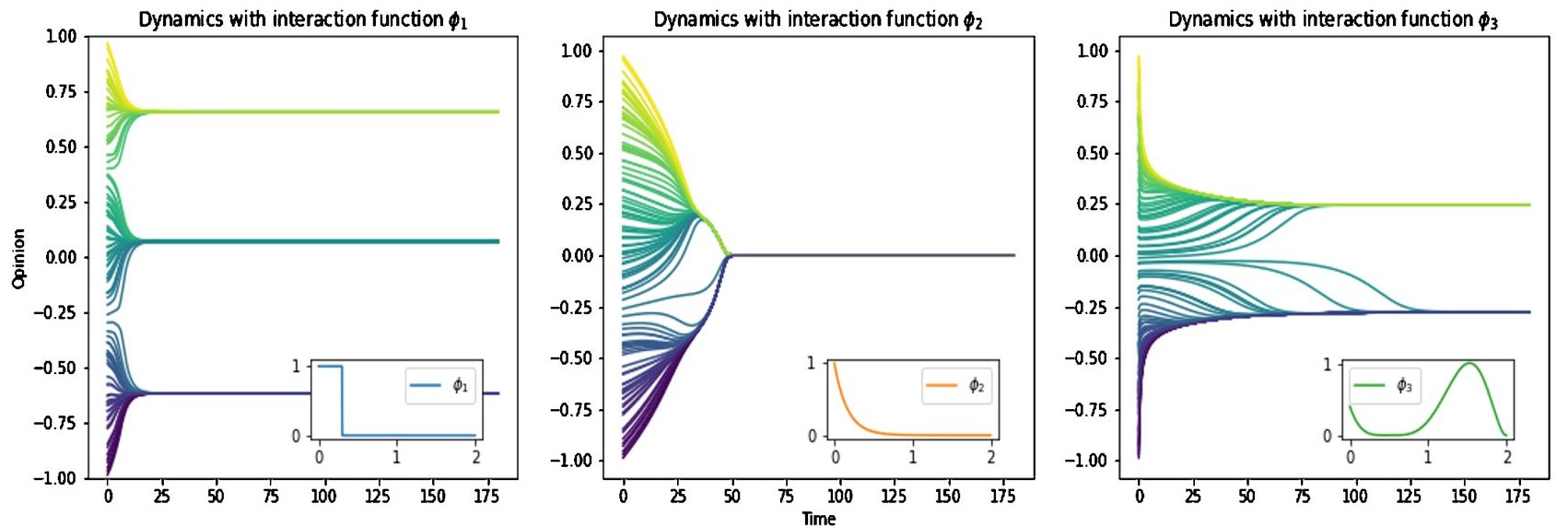}
    \caption{Examples of microscopic dynamics \eqref{eqn: HK-model} for different interaction functions $\phi$. Each interaction function is symmetric and takes values in $[0,1]$. We note different numbers of clusters appearing for each interaction function.}
    \label{fig:example uncontrolled miscroscopic system}
\end{figure}

We will first introduce microscopic dynamics for a finite population then derive the corresponding macroscopic dynamics. Consider a fixed population of $N$ agents with opinions $x_i(t) \in \mathcal{I}$. We assume that agents' opinions move closer when they interact, with the probability or strength of an interaction between individuals $i$ and $j$ being given by an \textit{interaction function} $\phi(x_j - x_i)$ \cite{nugent2024bridging}. Agents therefore effectively more towards a weighted average the opinions of those around them, using the interaction function as a weighting. This gives the following coupled ODE system, known as the Hegelsmann-Krause (HK) model \cite{hegselmann2015opinion},
\begin{equation}\label{eqn: HK-model}
    \dot{x}_i = \frac{1}{N}\sum_{j=1}^N\phi(x_j - x_i)(x_j - x_i) \,,
\end{equation}
for $i =1, ..., N$. Commonly, the interaction function is chosen such that $\phi(x_j - x_i)$ is symmetric and decreases with the distance $|x_j - x_i|$. A classical choice of interaction function, taken from the original HK model \cite{hegselmann2015opinion}, is the \textit{bounded confidence} interaction function. Here a confidence bound $R>0$ is specified and the interaction function is given by
\begin{align}\label{Eqn: bounded confidence kernel}
    \phi_R(x_j - x_i) = 
    \begin{cases}
        1 & \text{if } |x_j - x_i|<R \,, \\
        0 & \text{otherwise.}
    \end{cases}
\end{align}
This type of opinion formation model has been studied extensively, see for example \cite{motsch2014heterophilious, nugent2024bridging} and \cite{bernardo2024bounded} for a recent survey. Figure \ref{fig:example uncontrolled miscroscopic system} shows examples of the behaviour of our system \eqref{eqn: HK-model} for different interaction functions $\phi$. We see clustering when the interaction function is given by \eqref{Eqn: bounded confidence kernel}, corresponding to $\phi_1$ in Figure \ref{fig:example uncontrolled miscroscopic system}, and consensus when the interaction function is strictly positive, as is the case for $\phi_2$ in Figure \ref{fig:example uncontrolled miscroscopic system}. The final interaction function, $\phi_3$, considered in Figure \ref{fig:example uncontrolled miscroscopic system} is a so-called \textit{heterophilious} interaction function and assumes that while individuals communicate strongly with people who have opinions close to their own, they are also interested in people with differing opinions. This type of interaction is explored in depth in \cite{motsch2014heterophilious}.

We notice in Figure \ref{fig:example uncontrolled miscroscopic system} that the clustering of opinions is common in the HK model. In \cite{motsch2014heterophilious}, Motsch and Tadmor prove that given an interaction function $\phi$, the number of clusters can be quantified according to the multiplicity of the eigenvalue $\lambda=1$ of the adjacency matrix $a_{ij} = \frac{1}{N}\phi(x_j -x_i)$.

Models of populations can be either \textit{microscopic}, meaning they track each individual, or \textit{macroscopic}, meaning they track the population as a whole. The model introduced above is an example of a microscopic model. When considering a very large population size, a macroscopic \textit{mean-field} model is more appropriate (named for the velocity field generated by averaging over possible interactions). Instead of tracking the opinion of each agent individually we introduce an opinion density $\mu(x, t)$ which describes the density of the population with opinion $x\in\mathcal{I}$ at time $t\geq0$. The density is normalised such that
\begin{equation*}
    \int_{\mathcal{I}}\mu(x,t)\, dx = 1\,\quad \forall t\geq 0.
\end{equation*}
The time evolution of $\mu(x,t)$ can be derived either from considering the \textit{empirical measure} of the microscopic model in the limit as $N\to\infty$ \cite{canuto2012eulerian, goddard2022noisy} or from considering a \textit{Boltzmann-type description} of the system \cite{toscani2006kinetic,pareschi2013interacting}. For the purposes of this review, we will focus largely on the Boltzmann-type description as this has been more frequently employed in works on opinion manipulation. In \cite{toscani2006kinetic}, Toscani likens the random interactions between agents in a system to random interactions by gas molecules, the original application for the Boltzmann equation. We begin by considering binary interactions between two agents with opinions $x, x_*\in\mathcal{I}=[-1,1]$, producing opinions $x', x_*'$ given by
\begin{align}\label{eq:binary-interac-nocontrol}
\begin{split}
    x' &= x + \alpha \phi(x_*-x)(x_*-x) + \theta D(x),\\
    x_*' &= x_* + \alpha\phi(x - x_*)(x-x_*) + \theta_*D(x_*).
\end{split}
\end{align}
Here $\alpha\in(0,1/2)$ is a constant representing the extent to which individuals are willing to change opinion, $\theta,\theta_*$ are identically distributed random variables with mean 0 and variance $\sigma^2$ and $D(\cdot)$ is the local diffusion for a given opinion. We consider the time evolution of the density $\mu(x,t)$ in the weak formulation, i.e., when integrated against an arbitrary test function $\varphi(x)$. In the following equations, $\left<\cdot\right>$ denotes expectation with respect to random variables $\theta$ and $\theta_*$ and $\chi(\cdot)$ denotes the indicator function. The integro-differential equation for our system with binary interactions given by (\ref{eq:binary-interac-nocontrol}) is
\begin{equation}\label{Eqn: general Boltzmann}
    \frac{d}{dt}\int_{\mathcal{I}}\varphi(x)\mu(x,t)\, dx = (Q(\mu, \mu), \varphi),
\end{equation}
where
\begin{equation}
    (Q(\mu, \mu),\varphi) = \left<\int_{\mathcal{I}^2}\beta_{(x,x_*)\to(x', x'_*)}(\varphi(x') - \varphi(x))\mu(x,t)\mu(x_*,t)\, dx\, dx_*\right>
\end{equation}
and
\begin{equation*}
    \beta_{(x,x_*)\to(x',x_*')} = \chi\left(|x'|\in\mathcal{I}\right)\chi\left(|x_*'|\in\mathcal{I}\right)
\end{equation*}
is the transition rate from opinions $(x,x_*)$ to $(x',x_*')$. The transition rate ensures that the result of binary interactions (\ref{eq:binary-interac-nocontrol}) respect the bounds $x', x_*'\in\mathcal{I}$. In Toscani's derivation of the kinetic description of opinion dynamics models \cite{toscani2006kinetic}, a quasi-invariant limit is described where the frequency of interactions is taken to be very large while the impact of an interaction event is taken to be very small. Details of this limit can be found in Appendix \ref{app:quasi-invariant-boltz}. Employing this method, we can recover the following Fokker-Planck equation
\begin{equation}\label{eqn: mean-field-pde}
    \partial_t\mu(x,t) + \partial_x\left(\mathcal{P}[\mu]\mu(x,t)\right) = \frac{\tilde{\sigma}^2}{2} \partial_{xx}(D(x)^2\mu(x,t)),
\end{equation}
where $\mathcal{P}$ is an interaction force, describing the change in $\mu$ due to interactions taking place between agents, given by
\begin{equation*}
    \mathcal{P}[\mu](x) = \int_{\mathcal{I}}\phi(x_*-x)(x_*-x)\mu(x_*,t)\, dx_*.
\end{equation*}
In \eqref{eqn: mean-field-pde}, $\tilde{\sigma}^2$ is a rescaling of $\sigma^2$ under the quasi-invariant limit (see Appendix \ref{app:quasi-invariant-boltz} or \cite{toscani2006kinetic} for details). Note that this has much the same form as \eqref{eqn: HK-model} except with an integral over other opinions replacing the sum over other agents. The partial differential equation (\ref{eqn: mean-field-pde}) is complemented with no-flux boundary conditions (with the interpretation that agents cannot leave the opinion interval $\mathcal{I}$) and an initial condition $\mu(x,0) = \mu_0(x)$. Together with its boundary condition, equation \eqref{eqn: mean-field-pde} is often referred to as a mean-field description of opinion formation. 

Note that equation \eqref{eqn: HK-model} and correspondingly equation (\ref{eqn: mean-field-pde}) can be easily generalised to higher dimensions, interpreted as modelling multiple opinions evolving independently.

The practise of social manipulation often involves manipulating opinion formation such that a particular desired opinion is reached by members of a population. Mathematically, this problem can be formulated in the context of \textit{optimal control}. Optimal control is a common tool in engineering, for example in aerospace engineering where we might want to optimise a trajectory to minimise fuel usage but maximise range or ensure a safe landing. For introductory reading on optimal control, see \cite{kirk2004optimal, evans2005introduction, troltzsch2010optimal}.

In this context, we wish to optimise the evolution of individuals' opinions in order to reach some target opinion while minimising a cost associated with influencing agents, often interpreted as the cost of advertisement. The mechanism by which we influence individuals is a modelling choice and will be the focus of the rest of this section. In Section \ref{Section: Directly Affecting Opinions}, we will suppose that an external figure can directly affect opinions through advertising or media \cite{albi2014kinetic,albi2017mean} and in Section \ref{Section: Leaders and Followers} the population is divided into leaders and followers where only leaders can be affected by control \cite{albi2014boltzmann}. Section \ref{Sec: Liars} supposes that a member of the population is able to lie in a manner that is optimised to evolve the dynamics toward a given state \cite{glendinning2025what} and finally in Section \ref{Section: Network Control} we consider a formulation of control where individuals are placed on a dynamic network on which the control acts \cite{nugent2024steering}.

\subsection{Directly Affecting Opinions}\label{Section: Directly Affecting Opinions}

The first control approach we will consider is to directly affect opinions by applying a control $u_i$ to the dynamics of each agent's opinion. This approach is taken in \cite{albi2014kinetic, albi2017mean}. The purpose of applying a control is to try to move the final consensus that agents reach to some \textit{desired opinion}, $x_d\in\mathcal{I}$. In the microscopic model, this is written 
\begin{equation}\label{Eqn: directly affect dynamics}
    \dot{x}_i = \frac{1}{N}\sum_{j=1}^N\phi(x_j - x_i)(x_j - x_i) + u_i.
\end{equation}

Such a direct control is highly effective, as shown by the following Proposition \ref{Prop: Direct controllability}, whose proof can be found in \ref{app: proof of controllability}.
\begin{proposition} \label{Prop: Direct controllability}
    For any given initial conditions $x_i(0)$, any continuous interaction function $\phi:[-2,2]\rightarrow[0,1]$ and any target opinion $x_d \in [-1,1]$, there exists a control $u:\mathbb{R}^+ \rightarrow [-2,2]^{N \times N}$ such that \eqref{Eqn: directly affect dynamics} reaches consensus at $x_d$. Specifically this means that $x_i(t) \rightarrow x_d$ for all $i=1,\dots,N$ as $t\rightarrow\infty$. 
\end{proposition}

To encourage efficiency we may also choose a control that minimises the cost functional,
\begin{equation}\label{Eqn: directly affect cost functional}
    \mathcal{J}(x,u) = \int_0^T\left(\sum_{i=1}^N\frac{1}{2}(x_i - x_d)^2 + \frac{\nu}{2}\sum_{i=1}^N u_i^2\right)\, dt,
\end{equation}
where $x_d$ is the desired consensus opinion, $\nu>0$ is a regularisation parameter and $T>0$ is some fixed time horizon.

\begin{figure}[ht!]
    \centering
    \includegraphics[width=0.45\linewidth]{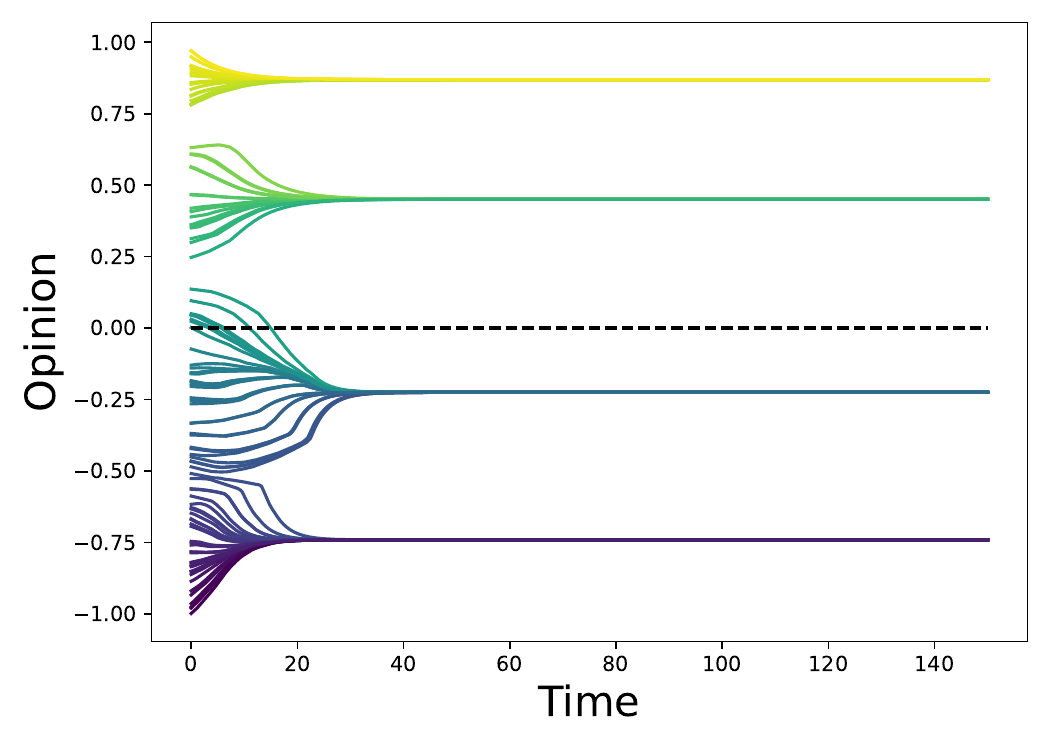}
    \includegraphics[width=0.45\linewidth]{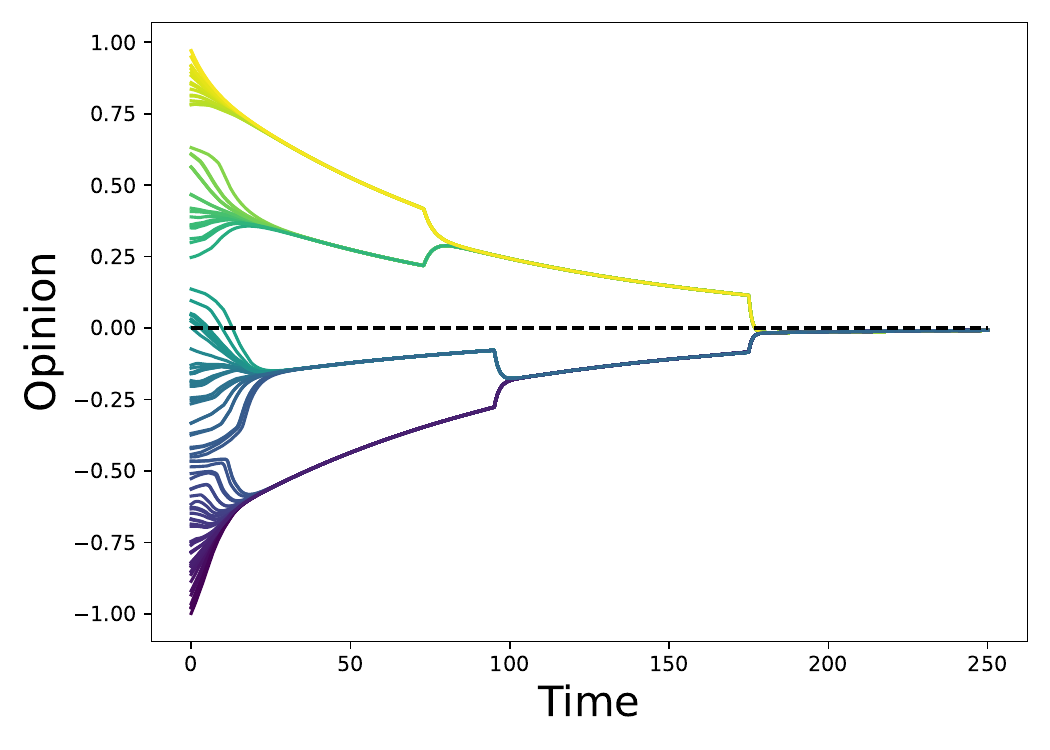}
    \caption{Example of an uncontrolled (left) and controlled (right) microscopic system with a bounded confidence interaction function with $R=0.2$. The target opinion, $x_d = 0$ is shown in the dashed black line. Here we consider $N = 100$ agents and the control is calculated according to cost function \eqref{Eqn: directly affect cost functional} with regularisation parameter $\nu = 1$.}
    \label{fig:ezample microscopic system directly affect}
\end{figure}

An example of a microscopic system of $N=100$ individuals being steered toward a goal opinion $x_d = 0$ by an instantaneous optimal control corresponding to cost functional \eqref{Eqn: directly affect cost functional} is shown in Figure \ref{fig:ezample microscopic system directly affect}. For details on the derivation of the instantaneous control for the microscopic model, see \cite{albi2014kinetic}. By comparing the left-hand figure, demonstrating the uncontrolled system, and the right-hand figure, showing the controlled system with $\nu =1$, we see that controlling opinions by directly affecting their evolution is a highly effective strategy for achieving some goal opinion $x_d$.

\begin{remark}
A derivation of control directly from a mean-field model is performed in \cite{albi2017mean} using established methods for optimal control of partial differential equations \cite{troltzsch2010optimal}. The cost functional $\mathcal{J}(x,u)$ in equation \eqref{Eqn: directly affect cost functional} becomes
$$\tilde{\mathcal{J}}(\mu,f) = \int_0^T\left(\frac{1}{2}\int(x-x_d)^2\mu(x,t)\, dx + \frac{\nu}{2}\int (f(x,t))^2\mu(x,t)\, dx\right)\, dt$$
where $f$ is an infinite dimensional generalisation of the control $u_i$ discussed in equation \eqref{Eqn: directly affect dynamics} and $\mu(x,t)$ is an opinion density satisfying
\begin{equation*}
    \partial_t\mu + \partial_x\left((\mathcal{P}[\mu]+ f)\mu\right) = \frac{\sigma^2}{2}\partial_{xx}\mu.
\end{equation*}
Existence of an optimal control $f$ is proven in \cite{albi2017mean}.
\end{remark}

Analysis of the success of our control is carried out by considering the large population limit, $N\to \infty$. This limit is carried out by considering a Boltzmann-type equation for the evolution of opinions, similar to that in Section \ref{Section: Introduction continuous models}. We first introduce an opinion density for the system $\mu(x,t)$, normalised such that
\begin{equation*}
    \int_{\mathcal{I}}\mu(x,t)=1, \quad \forall t\geq 0.
\end{equation*}
Under the conditions of Proposition 3.1 in \cite{albi2014kinetic}, we can express the time evolution of $\mu(x,t)$ via the Boltzmann equation \eqref{Eqn: general Boltzmann} with the interaction integral
\begin{equation}\label{Eqn: general interaction integral control}
    (Q(\mu, \mu), \varphi) = \eta\left<\int_{\mathcal{I}^2}(\varphi(x') - \varphi(x))\mu(x,t)\mu(x_*,t)\, dx\, dx_*\right>.
\end{equation}
Here, $x'$ is the result of a binary interaction between $x$ and $x_*$, similarly to \eqref{eq:binary-interac-nocontrol} but with the addition of the control term $u$ which can be explicitly calculated as an instantaneous control, detailed in \cite{albi2017mean}. Additionally, $\eta>0$ is a constant rate and $\left<\cdot\right>$ denotes expectation.

After performing a quasi-invariant limit, similar to that explained in Appendix \ref{app:quasi-invariant-boltz}, the authors obtain the following Fokker-Planck equation for the evolution of the density $\mu(x,t)$,
\begin{equation}\label{eqn: fokker-planck-direct-affect}
    \partial_t\mu(x,t) + \partial_x(\mathcal{P}[\mu]\mu(x,t)) + \partial_x(\mathcal{H}[\mu]\mu(x,t)) = \frac{\tilde{\sigma}^2}{2}\partial_{xx}\left(D(x)^2\mu(x,t)\right),
\end{equation}
where $D(x)$ is a diffusion function, $\tilde{\sigma}^2$ is a rescaling of the variance of the noise applied to binary interactions between agents. Further,
\begin{align}\label{Eqn: common interaction drift}
    \mathcal{P}[\mu](x) &= \int_{\mathcal{I}}\phi(x_* - x)(x_* - x)\mu(x_*,t)\, dx_*, \\
    \mathcal{H}[\mu](x) &= \frac{1}{\kappa}(x_d -x),
\end{align}
where $\kappa$ is a rescaling of the regularisation parameter $\nu$. We recognise $\mathcal{P}[\mu](x)$ as the drift term toward the mean of the neighbours of an agent with opinion $x$, weighted by interaction function $\phi(x_*-x)$, discussed in Section \ref{Section: Introduction continuous models}. The term $\mathcal{H}[\mu](x)$ represents the influence of the control, which is a drift toward $x_d$.

In cases where the interaction function $\phi$ and diffusion $D(x)$ are very simple functions, it is possible to analyse the steady state of this equation. In these cases, it can be shown that the mean of the system tends to $x_d$ with a sharper peak in the distribution $\mu$ for small values of the rescaled regularisation parameter $\kappa$ \cite{albi2014kinetic}.

In general, equation (\ref{eqn: fokker-planck-direct-affect}) is analytically intractable so the authors resort to numeric methods. For example, the Boltzmann Monte Carlo method simulates binary interactions in order to recover an approximation of the density profile $\mu(x,t)$ \cite{pareschi2013interacting}, proven to converge to the true distribution $\mu(x,t)$ as the number of particles $N\to\infty$ \cite{wagner1992convergence}. It is shown that the control works most effectively for small values of the regularisation parameter $\kappa$ and that for extreme values of $x_d$, for example $x_d = 0.95$, the control is less effective \cite{albi2014kinetic}, meaning that the time taken for the system to converge to a goal opinion $x_d$ is substantially greater.

The results of Albi et al. show that directly affecting opinions is an effective method for achieving consensus at a desired opinion. In this context, the control $u_i$ must be interpreted as an external forcing which could come in the form of social media misinformation or advertisement. Various details of the model have been explored in \cite{albi2017mean}, including a derivation of optimality conditions for controls beginning from the PDE model \eqref{eqn: mean-field-pde}. 

Controlling a population by directly affecting opinions is a very flexible and effective strategy for manipulation. In personalising the control for each agent $i$, which can be thought of as targeted advertisement or algorithm personalisation, we can steer each individual in the population toward our desired opinion without having a large effect on the trajectory of nearby individuals. Indeed, in the microscopic case \eqref{Eqn: directly affect dynamics}, the instantaneous optimal control takes into account the influence that each agent has on one another. However, in real life, it may not be feasible or one may not have the resources to tackle agents individually in the manner described above. In this case, it seems more reasonable to assign a control strategy that targets sections of the opinion space so the same control will be assigned to any agent with an opinion within a certain subinterval of $\mathcal{I}$. 

\subsection{Leaders and Followers}\label{Section: Leaders and Followers}

In this section, we consider a model where some agents act with a manipulative strategy. A common feature of opinion formation models is the separation of the population into \textit{followers} and \textit{leaders} \cite{albi2014boltzmann, during2015opinion}. The opinion of a follower is influenced by the opinions of other followers and the opinions of leaders. On the other hand, the opinion of a leader is influenced only by the other leaders and here an additional control variable $u$. We suppose that our system contains $N_F$ following agents and $N_L$ leaders. The followers have opinions $x_i\in\mathcal{I}$ for $i=1, ..., N_F$ and leaders have opinions $\tilde{x}_k\in\mathcal{I}$ for $k = 1, ..., N_L$. Opinions of leaders and followers evolve with respect to the following sets of ODEs,
\begin{subequations}\label{eqn: full leaders and followers microscopic}
\begin{align}
    \dot{x}_i &= \frac{1}{N_F}\sum_{j=1}^{N_F}\phi(x_j - x_i)(x_j - x_i) + \frac{1}{N_L}\sum_{k=1}^{N_L}\phi(\tilde{x}_k - x_i)(\tilde{x}_k - x_i), & x_i(0) = x_{i,0} \label{eqn: LF_microscopic_followers} \\
    \dot{\tilde{x}}_k &= \frac{1}{N_L}\sum_{j=1}^{N_L}\phi(\tilde{x}_j- \tilde{x}_k)(\tilde{x}_j - \tilde{x}_k) + u, & \tilde{x}_k(0) = \tilde{x}_{k,0}. \label{eqn: LF_microscopic_leaders}
\end{align}
\end{subequations}
The control term $u$ characterises the strategy of the leaders and is determined by solving the following optimal control problem:
$$ u = \arg\min_{u\in\mathbb{R}}\{\mathcal{J}(u, \mathbf{x},\mathbf{\tilde{x}})\},$$
where
\begin{equation}\label{eqn:cost_functional_LF}
    \mathcal{J}(u, \mathbf{x}, \mathbf{\tilde{x}}) = \frac{1}{2}\int_0^T\left(\frac{\psi}{N_L}\sum_{k=1}^{N_L}(\tilde{x}_k-x_d)^2 + \frac{\gamma}{N_L}\sum_{k=1}^{N_L}(\tilde{x}_k - m_F)^2\right)\, ds + \int_0^T\frac{\nu}{2}u^2\, ds.
\end{equation}
Here, $\mathbf{x}$ and $\mathbf{\tilde{x}}$ are vectors representing the opinions of followers and leaders respectively, $x_d$ is a desired opinion, $m_F = \frac{1}{N_F}\sum_{j=1}^{N_F}x_j$ is the mean opinion of the followers and $T$ represents the final time horizon. The parameters $\psi, \gamma$ and $\nu$ are all non-negative constants. We impose that $\psi + \gamma =1$ to balance the considerations of the control to move the leaders' opinions toward the desired state and to remain close to the mean opinion of the followers. As usual, $\nu\geq 0$ is a regularisation constant for the control, penalising the action of the control $u$ for large values of $\nu$.

Similarly to Section \ref{Section: Directly Affecting Opinions}, the natural extension of this model is to consider the large population limit, where $N_F, N_L\to \infty$ \cite{albi2014boltzmann}. We first introduce separate densities for the leaders and the followers. On some interval $\mathcal{I}$, we introduce opinion variables $x, \tilde{x}\in\mathcal{I}$ and corresponding densities $\mu_F(x,t)$ and $\mu_L(\tilde{x},t)$ representing the probability density of followers with opinion $x$ at time $t\geq0$ and leaders with opinion $\tilde{x}$ respectively. These densities are normalised such that
$$\int_{\mathcal{I}} \mu_F(x,t)\, dx = 1, \quad \int_{\mathcal{I}} \mu_L(x,t)\, dx = \rho,  \quad \forall t\geq 0,$$
where $\rho\geq 0$ is a constant parameter. The parameterisation of the volume of leaders, $\rho$ in $\mathcal{I}$ allows us to consider scenarios where the relative population of leaders to followers is small. We then have two resulting Boltzmann equations, one for the followers and one for the leaders:
\begin{equation}\label{Eqn: boltzmann leaders followers}
    \begin{split}
    \frac{d}{dt}\int_{\mathcal{I}}\varphi(x)\mu_F(x,t)\, dx &= (Q_{F}(\mu_F, \mu_F), \varphi) + (Q_{FL}(\mu_F, \mu_L), \varphi),\\
    \frac{d}{dt}\int_{\mathcal{I}}\varphi(\tilde{x})\mu_L(\tilde{x},t)\, d\tilde{x} &= (Q_{LL}(\mu_L, \mu_L), \varphi).
    \end{split}
\end{equation}
Here, $(Q_{F}(\mu_F, \mu_F), \varphi)$ encodes the change in $\mu_F$ due to interactions amongst followers, $(Q_{FL}(\mu_F, \mu_L), \varphi)$ encodes the change in $\mu_F$ due to interaction of followers with leaders and $(Q_{L}(\mu_L, \mu_L), \varphi)$ encodes the change in $\mu_L$ due to interactions amongst leaders, taking into account the control. Under the assumptions outlined in Proposition 3.1 in \cite{albi2014boltzmann}, we have that $(Q_F(\mu_F,\mu_F),\varphi)$, $(Q_{FL}(\mu_F, \mu_L),\varphi)$ and $(Q_L(\mu_L, \mu_L), \varphi)$ all take the same form as equation \eqref{Eqn: general Boltzmann} with the respective binary interactions between each population. Similarly to in Section \ref{Section: Directly Affecting Opinions}, we introduce constant rates $\eta_F, \eta_{FL}, \eta_L>0$ for each interaction integral. We can take the quasi-invariant limit of interactions, detailed in \cite{albi2014boltzmann}, resulting in the coupled Fokker-Planck equations
\begin{align*}
    \partial_t\mu_F(x,t) + \partial_x\left(\left(\frac{1}{c_F}\mathcal{P}_F[\mu_F](x) + \frac{1}{c_{FL}}\mathcal{P}_{FL}[\mu_L](x)\right)\mu_F(x)\right) &=\frac{\tilde{\sigma}^2}{2}\left(\frac{1}{c_T} + \frac{\rho}{c_{FL}}\right)\partial_{xx}\left( D^2(x)\mu_F(x)\right),\\
    \partial_t\mu_L(\tilde{x},t) + \partial_{\tilde{x}}\left(\left(\frac{\rho}{c_L}\mathcal{H}[\mu_L](\tilde{x}) + \frac{1}{c_L}\mathcal{P}_L[\mu_L](\tilde{x})\right)\mu_L(\tilde{x},t)\right) &= \frac{\tilde{\sigma}^2}{2}\frac{\rho}{c_L}\partial_{\tilde{x}\tilde{x}}\left(D^2(\tilde{x})\mu_L(\tilde{x}, t)\right).
\end{align*}
Here, $\frac{1}{c_F}, \frac{1}{c_{FL}}, \frac{1}{c_L}$ are rescalings of the rate parameters $\eta_F, \eta_{FL}$ and $\eta_L$ and $\tilde{\sigma}^2$ is a rescaling of the variance of the noise applied to binary interactions between agents. Furthermore,
\begin{align*}
    \mathcal{P}_{(\cdot)}[\mu_{(\cdot)}](x) &= \int_{\mathcal{I}}\phi(x_* - x)(x_* - x)\mu_{(\cdot)}(x_*,t)\, dx_*,\\
    \mathcal{H}[f_L](\tilde{x}) &= \frac{2\psi}{\kappa}(\tilde{x}+m_L(t) - 2x_d) + \frac{2\gamma}{\kappa}(\tilde{x}+m_L(t) -2m_F(t)),
\end{align*}
where $\kappa>0$ is a rescaling of the regularisation parameter $\nu$. The terms $\mathcal{P}_{F}$ and $\mathcal{P}_{FL}$ in the convolution represent drifts by the followers toward the mean opinion of the followers and the leaders, respectively. The term $\mathcal{P}_L$ represents a drift in the opinion density of the leaders toward the mean opinion of the leaders. Finally, the term $\mathcal{H}[f_L](\tilde{x})$ represents the action of the control on the leaders' opinions. The first term in $\mathcal{H}[f_L](\tilde{x})$ is a shift toward the goal opinion $x_d$, and the second is a shift toward the mean opinion of the followers, $m_F(t)$.

In simple cases, steady states of the Fokker-Planck equations can be analysed analytically, with results presented in \cite{albi2014boltzmann}. Indeed, an analysis of the Boltzmann equation under the quasi-invariant limit, seen in Section 4 of \cite{albi2014boltzmann}, shows that the steady state of the dynamics of the mean opinions for each population is the goal opinion $x_d$ in the simple case where $\phi(x' - x) = 1$ for all $x',x\in\mathcal{I}$. Similarly to when opinions are directly affected in Section \ref{Section: Directly Affecting Opinions}, the Boltzmann Monte Carlo method is a useful tool for obtaining numeric results \cite{pareschi2013interacting, albi2014boltzmann}. 

Simple extensions to the model also yield interesting results, such as the long term distribution of followers in the presence of multiple sets of leaders. Authors in \cite{albi2014boltzmann} also explore time-dependent strategies, where the constants $\psi$ and $\gamma$ in the cost functional (\ref{eqn:cost_functional_LF}) are replaced by functions $\psi(t)$ and $\gamma(t)$ varying in time. This allows more control over when a radical or a populistic strategy for the control is most helpful.

The concept of leaders in a population is often combined with that of assertiveness or stubbornness. In \cite{during2015opinion}, authors take the assertiveness of agents to be a variable and call individuals with high levels of assertiveness leaders. Despite no control being placed on this model and no desired opinion stated, it is observed that allowing agents to have assertiveness of their opinions enhances opinion clustering within a population. Similarly, in \cite{during2009boltzmann}, stubborn leaders with high assertiveness are treated as a subset of the population with a resulting Boltzmann description similar to equation \eqref{Eqn: boltzmann leaders followers}. Again, no control is placed on the system but in choosing an initial distribution of leaders' opinions we can steer followers opinions toward a target distribution. Finally, in \cite{during2024breaking}, authors steer a leaders and followers formulation away from consensus by introducing control on a network which in the mean-field limit becomes a graphon. We will see further discussion of network control on the microscopic level in Section \ref{Section: Network Control}.

A leaders and followers formulation is a realistic formulation for consensus control. The assumption that only certain members of the population can be influenced but can steer the population is similar to what is seen in real life with social media influencers who are paid to promote products (analogues to our control) which leads to those products becoming more popular in a wider population. Mathematically, the leaders and followers formulation is quite similar to the idea of directly affecting opinions, explored in Section \ref{Section: Directly Affecting Opinions}, however there is less freedom in that followers cannot be directly affected by the control and the assumptions made in \eqref{eqn: full leaders and followers microscopic} mean that the same control is used to influence all leaders which is less effective than individually targeting each leader.

\subsection{Control through Lying}\label{Sec: Liars}

Another method for adding control to the HK model \eqref{eqn: HK-model} is to suppose that there exists an agent (or group of agents) who can present \textit{apparent opinions} or \textit{lies} to each member of the population. This section follows the work of Glendinning et. al in \cite{glendinning2025what}. The lies are chosen such that consensus of the truth-telling population at a desired opinion is achieved in as short an amount of time as possible, with some regularisation constraint. Under the assumption that one agent, labelled agent 1 without loss of generality, can lie by presenting an apparent opinion $u_i$ to a truth-telling agent with label $i$, the microscopic formulation becomes
\begin{equation}
    \dot{x}_i = \frac{1}{N}\left(\phi(u_i - x_i)(u_i-x_i)+ \sum_{j=2}^N\phi(x_j-x_i)(x_j-x_i)\right)
\end{equation}
for $i=2, ..., N$ and with initial conditions $x_i(0) = x_i^0$. We will additionally assume that the lying agent does not change their true opinion in time and that the liar's initial opinion is their goal opinion, so
\begin{equation*}
    \dot{x}_1 = 0, \quad x_1(0) = x_d.
\end{equation*}
In a similar manner as in the previous sections, we will impose a cost on the control, written as the deviation of the liar from his/her true opinion. The lies $u_i$ are chosen to minimise
\begin{equation} \label{Eqn: Lying cost function}
    \mathcal{J}(\mathbf{x},\mathbf{u}) = \int_0^T \left(\sum_{i=1}^N\frac{1}{2}(x_i-x_d)^2 + \frac{\nu}{2}\sum_{i=2}^N(u_i-x_d)^2\right)\, ds,
\end{equation}
with the restriction that $u_i\in\mathcal{I}$ for each $i=2, ..., N$. This ensures that the liar can only express opinions inside the opinion interval $\mathcal{I}$. Here, $T>0$ is an arbitrary time horizon and $\nu>0$ is a constant regularisation parameter. In \cite{glendinning2025what}, Glendinning et. al explore multiple regularisation strategies, inspired by social constraints. For example, it might be intuitive to penalise the liar for not being consistent in the lies they tell to different agents. The resulting instantaneous control from each regularisation is compared for simple cases of interaction function $\phi$.

Again, it is useful to consider a Boltzmann-type description of the system. Similarly to the leaders and followers formulation described in Section \ref{Section: Leaders and Followers}, we need to consider the liar and the truth-tellers as two separate populations. Let $\mu_T(x,t)$ denote the density of truth-tellers with opinion $x$ at time $t\geq0$. We normalise this density to 1 as in the previous sections. We also introduce a density for the liar's opinion $\mu_L(x,t)$, normalised to a constant parameter $0<\rho\leq1$ representing the relative influence of the liar on a particular agent compared to the influence of the general population of truth-tellers on said agent. Of course, since the liar's opinion does not change in time we have
$$ \mu_L(x,t) = \rho\delta(x-x_d),\quad \forall t\geq0.$$
The Boltzmann equation for the truth-telling population is then given by
\begin{equation}
    \frac{d}{dt}\int_{\mathcal{I}}\varphi(x)\mu_T(x,t)\, dx = (Q_{T}(\mu_T,\mu_T),\varphi) + (Q_L(\mu_T, \mu_L), \varphi),
\end{equation}
where interaction integrals $(Q_T(\mu_T,\mu_T),\varphi)$ and $(Q_L(\mu_T, \mu_L),\varphi)$ represent interactions between truth-tellers and interactions between the truth-tellers and the liar respectively. Under the assumptions outlined in Proposition 4.2 in \cite{glendinning2025what}, we have that $(Q_T(\mu_T, \mu_T),\varphi)$ and $(Q_L(\mu_T,\mu_L),\varphi)$ take the form of equation \eqref{Eqn: general Boltzmann} with respective binary interactions corresponding to interactions among truth-tellers and interactions between the truth-tellers and the liar. We also introduce positive rates $\eta_T$ and $\eta_L$ for each interaction integral. Under a quasi-invariant limit, we recover the following Fokker-Planck equation for the density of truth-telling opinions
\begin{equation}\label{Eqn: Fokker-Planck Liars}
    \partial_t\mu_T(x,t) + \partial_x(\mathcal{P}[\mu_T](x)\mu_T(x,t)) + \partial_x(\mathcal{H}[\mu_T](x)\mu_T(x,t)) = \frac{\tilde{\sigma}^2}{2}\left(\frac{1}{c_T}+\frac{\rho}{c_L}\right)\partial_{xx}(D(x)^2\mu_T(x,t)).
\end{equation}
Here, $1/c_T$ and $1/c_L$ are rescaled versions of the interaction rates $\eta_T$ and $\eta_L$, $\tilde{\sigma}^2$ is a rescaling of the noise applied to binary interactions, $\mathcal{P}[\mu_T](x)$ is given by \eqref{Eqn: common interaction drift} and
$$ \mathcal{H}[\mu_T](x) = \phi(u(x) - x)(u(x)-x)$$
where $u(x)$ lie told to an individual with opinion $x$. As in previous formulations, $\mathcal{P}[\mu_T](x)$ represents the influence of the truth-telling population on a truth-teller with opinion $x$, whereas $\mathcal{H}[\mu_T](x)$ represents the influence of the liar. Taking $u(x)$ to be the quasi-invariant limit of the instantaneous binary control, we have
\begin{equation}\label{Eqn: instantaneous quasi-invariant control}
    u(x) = x_d - \frac{1}{\kappa}(x-x_d)\partial_u(\phi(u(x)-x)(u(x)-x)).
\end{equation}
Here $\kappa$ is a rescaling of the regularisation parameter $\nu$. Note that this is the only formulation we have considered so far where the instantaneous control depends on a derivative of the interaction function. In simple cases, such as $\phi = \phi(x)$, where $x$ is the opinion of a truth-telling agent, equation \eqref{Eqn: instantaneous quasi-invariant control} is analytically tractable and steady state solutions can be derived.

\begin{figure}[ht!]
    \centering
    \includegraphics[width=0.48\linewidth]{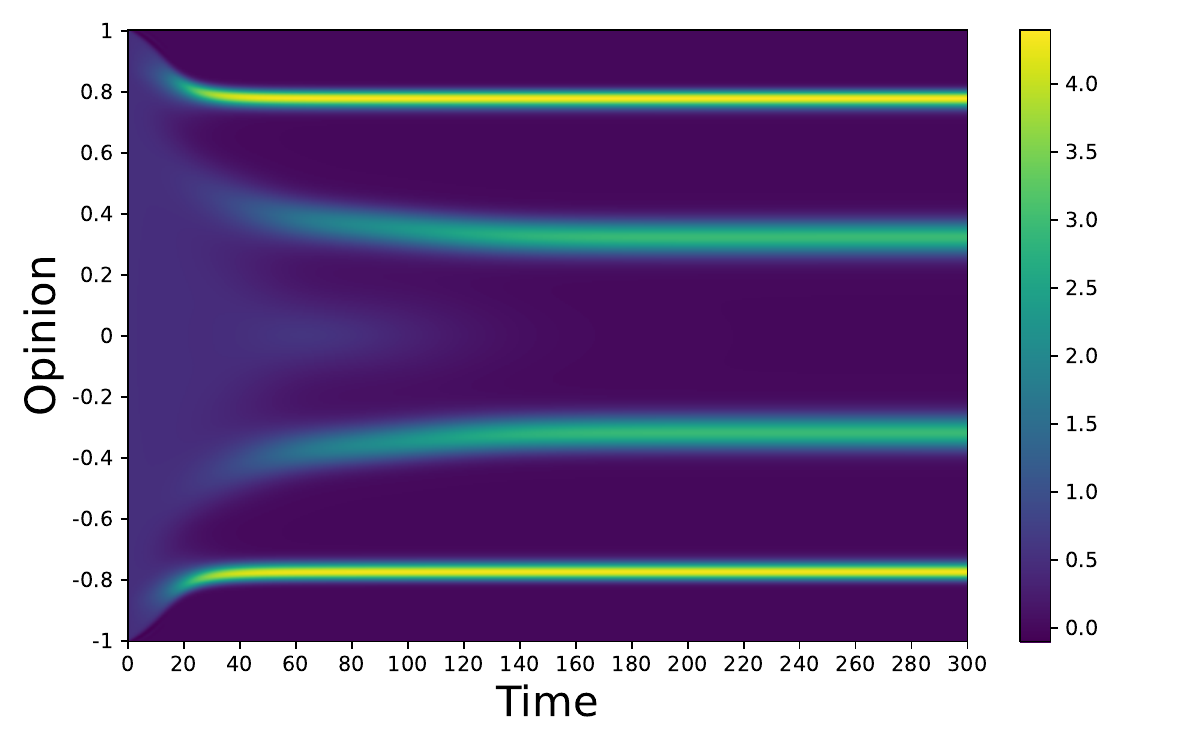}
    \includegraphics[width=0.48\linewidth]{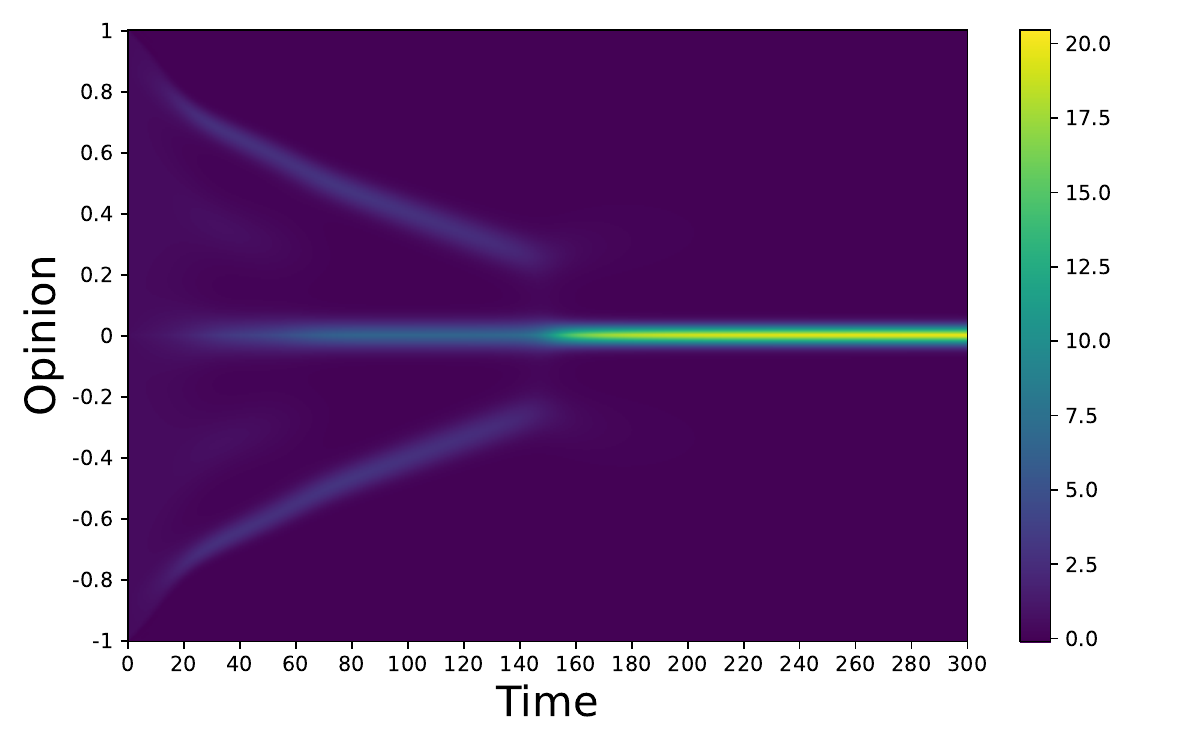}
    \caption{Uncontrolled (left) and controlled (right) solutions to the Fokker-Planck equation \eqref{Eqn: Fokker-Planck Liars} in the case where the liar's goal opinion $x_d = 0$. In the left-most plot, we have $\rho=0$ and the population interacts without any external influence and we observe clustering. In the right-most plot, we take $\rho=0.1$ and $\kappa = 0.01$ and observe that the population is brought to consensus at the liar's desired opinion. In both plots, the interaction function is equal to 1 for $|x-x'|<r_1 = 0.1$ and equal to 0 for $|x-x'|>r_2=0.2$ with smooth interpolation between these two values. The right-most plot is reproduced with permission from \cite{glendinning2025what}.}
    \label{fig:example macroscopic bounded con liars}
\end{figure}

In \cite{glendinning2025what}, Glendinning et. al develop a numeric method for discretising $u(x)$ in the case where $\phi(x_* - x)$ is a smoothed version of the bounded confidence kernel \eqref{Eqn: bounded confidence kernel}. It is found that for sufficiently small $\kappa$ and any $\rho>0$, it is possible to bring the entire truth-telling density to be concentrated around $x_d$. Figure \ref{fig:example macroscopic bounded con liars} shows an example of the uncontrolled macroscopic system (\eqref{Eqn: Fokker-Planck Liars} with $\rho=0$) and an example of the system where the liar is active. We observe that for a small value of $\kappa$, the liar is able to communicate with the entire population and achieves consensus at $x_d=0$. In contrast, for the uncontrolled system we observe opinion clusters and a very low density of individuals with opinion $x_d$ at large times.

The authors also discuss the case where there are multiple lying populations, each with a different goal opinion. We define a weak liar to be a liar with large regularisation parameter $\kappa$, and a strong liar to have small $\kappa$. In this case, an interesting feature of the model is that weak liars perform better in the presence of a strong liar with a different goal opinion. Indeed, individuals that would previously be unreachable by the weak liar can be influenced by the strong liar into the confidence region of the weak liar. Competition and coordination between liars can be seen as analogous to game theory approaches to manipulation in opinion dynamics which will be discussed in Section \ref{sec: coord and coop}.

The formulation of optimal control through lying is a highly socially intuitive mechanism for manipulating a population. Altering the model by penalising the liar for deviating too far from their opinion, telling too many different lies to too many people or changing their apparent opinion too much in time leads to a rich array of dynamics and optimal strategies of the liar. In contrast to the formulations presented in Sections \ref{Section: Directly Affecting Opinions} and \ref{Section: Leaders and Followers}, the population is not being controlled by any external figure but rather by an agent within the population, meaning that these results could be verified experimentally.

\subsection{Network Control} \label{Section: Network Control}

The final opinion manipulation mechanism we consider is that of network control. In this case opinion manipulation is not direct, there is no external or leader influence on individuals' opinions, instead it is their underlying social network that is affected. Such a control could be enacted in the real world by, for example, manipulating the popularity or viewership of social media posts. 

It is already common to include a social network over which individuals interact in models of opinion formation \cite{amblard2004role,gabbay2007effects}, and increasingly common in both opinion dynamics and network science more broadly to study dynamic networks \cite{nugent2023evolving,berner2023adaptive,kozma2008consensus,iniguez2009opinion}. 

Before introducing the model considered in this section we first establish some notation used throughout for networks. A network (or graph) $G = (V,E)$ is comprised of a set of nodes ($V$) and edges ($E$). In all the works considered here the node set will be the set of agents $V = \{1,\dots,N\}$. The edge set is a list of pairs $(i,j)$ indicating the existence of a connection from $i$ to $j$. We typically consider directed networks, meaning that an edge from $i$ to $j$ does not necessarily imply the existence of an edge from $j$ to $i$. Moreover, edges are often weighted with the weight of the edge from $i$ to $j$ denoted $w_{ij}$. These edge weights form the weighted adjacency matrix $W$, with $w_{ij}=0$ if there is no edge from $i$ to $j$. This adjacency matrix gives an alternative description of the network, and the term network is used to refer to both the graph $G=(V,E)$ and its corresponding adjacency matrix $W$. 

In \cite{nugent2023evolving}, Nugent et al. extend the HK model \eqref{eqn: HK-model} to include a dynamic network, later developing this into a control problem in \cite{nugent2024steering}. As with \eqref{eqn: HK-model}, the model begins with $N$ ordinary differential equations for individuals' opinions $x_i \in [-1,1]$. In addition to this we now include $N^2$ ordinary differential equations for the network edge weights $w_{ij}\in[0,1]$. The dynamics are, given by 
\begin{subequations} \label{eqn: ODE system} 
 \begin{align} 
    \frac{dx_i}{dt} &= \frac{1}{d_i} \sum_{j=1}^N w_{ij}\, \phi(|x_j - x_i|) \,(x_j - x_i) \quad & i=1,\dots,N \,\label{eqn: opinion ODE}, \\
    \frac{dw_{ij}}{dt} &= s(u_{ij}) \, \big( \ell(u_{ij}) - w_{ij} \big) \quad & i,j=1,\dots,N \,, i\neq j \,. \label{eqn: weight ODE}
\end{align}
\end{subequations}
Now $u_{ij}\in[-1,1]$ is a control variable affecting the evolution of the edge weight $w_{ij}$. The functions $s:[-1,1]\rightarrow[0,\mathcal{S}]$ and $\ell:[-1,1]\rightarrow[0,1]$ control the speed and direction of controls respectively, and are given by 
\begin{align*}
    s(u) &= \mathcal{S} u^2 \,,\\
    \ell(u) &= \frac{1}{2}(u+1) \,,
\end{align*}
where $\mathcal{S}>0$ is a constant maximum control speed. Hence, setting a control to $u_{ij}=-1$ reduces the edge weight $w_{ij}$ towards zero, setting a control to $u_{ij} = 0$ has no effect, setting a control to $u_{ij}=1$ increases an edge weight towards one, with intermediate control values targetting intermediate edge weights. 

This setup allows the control to have a significant impact on the structure of the network, but the choice of $\mathcal{S}$ limits the speed at which this effect can take place. This is somewhat different to other network controls which instead instantaneously switch edge on or off, see for example \cite{castiglioni_election_2020}, an approach is more closely related to that of ‘edge-based’ controls or a ‘decentralised adaptive strategy’ which has been applied to other interacting particle systems \cite{de2009decentralized,delellis2010synchronization,rajapakse2011dynamics,yu2012distributed}. This approach prescribes edge weight dynamics that are designed to encourage consensus, similar to those introduced to encourage cooperation in Section \ref{sec: coord and coop}, but crucially does not introduce a control variable that determines edge weight dynamics. A similar approach to \cite{nugent2023evolving} was taken by Piccoli and Duteil in \cite{piccoli2021control}, where they introduce a control on individuals `masses', denoted $m_i$ which determine their reputation. This effectively acts as a directed network with the specific structure that $w_{ij} = m_j$ for all $i$.

The goal considered in both \cite{piccoli2021control} and \cite{nugent2024steering} is to bring the entire population to consensus at a target opinion $x_d$. It is shown in Theorem 3.1 of \cite{nugent2024steering} that such a control is possible for \eqref{eqn: ODE system} under mild assumptions on the interaction function and initial conditions, provided that the control acts sufficiently quickly. However, the most effective controls are obtained through the study of optimal controls, aiming to minimise the following cost functional: 
\begin{equation} \label{Eqn: network control cost functional}
    \mathcal{J}(u) = \int_0^T \alpha \sum_{i=1}^N \sum_{j=1}^N u_{ij}(s)^2 + \beta \sum_{i=1}^N \big(x_i(s) - x_d \big)^2 \,\,ds \,. 
\end{equation}
The first term describes the cost of control while the second term describes the total distance from the target opinion $x_d$. Integrating both these costs over time encourages the application of minimal controls while bringing the population towards the target as quickly as possible, with the relative weights $\alpha$ and $\beta$ balancing these two objectives. In this case the optimal controls are bang-bang, meaning they are either zero or take the maximum value of $\pm1$. 

\begin{figure}[ht!]
    \centering
    \includegraphics[width=0.8\linewidth, trim = {2cm 2cm 2cm 2cm}, clip]{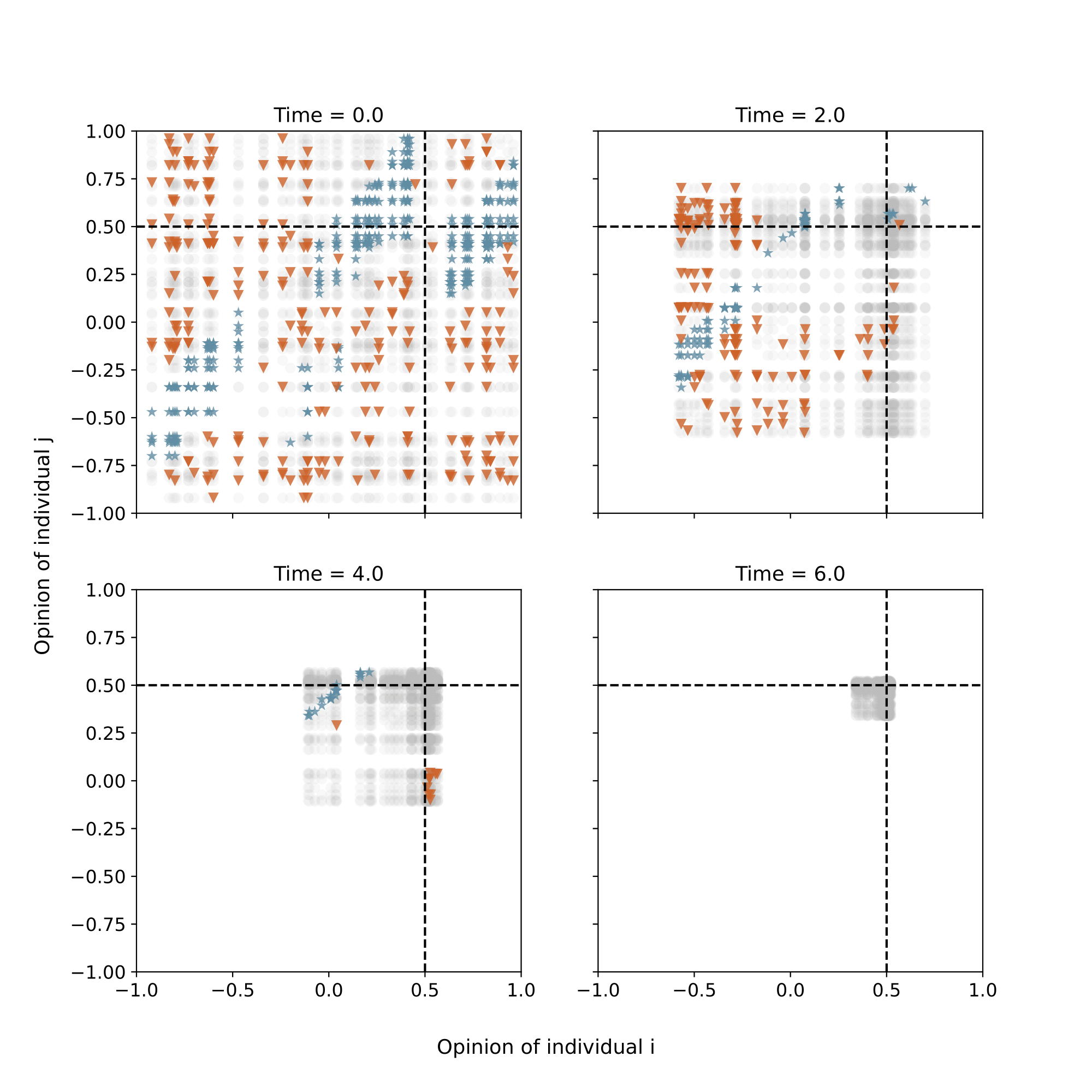}
    \caption{Reproduced with permission from \cite{nugent2024steering}. Example optimal control of \eqref{eqn: ODE system} with cost function \eqref{Eqn: network control cost functional}. At each timepoint a blue star, grey circle or red triangle is placed at $(x_i,x_j)$ is the control $u_{ij}$ is equal to $+1,0$ and $-1$ respectively. The target opinion $x_d=0.5$ is shown with dashed lines. The control is successful as all points are converging towards this target.}
    \label{fig:example optimal network controls}
\end{figure}

Figure \ref{fig:example optimal network controls} shows an example of optimal network controls obtained for this problem using a bounded confidence interaction function. At each timepoint a blue star, grey circle or red triangle is placed at $(x_i,x_j)$ is the control $u_{ij}$ is equal to $+1,0$ and $-1$ respectively. The target opinion $x_d=0.5$ is shown with dashed lines. We can see that the control is successful as all points are converging towards this target.

At times $t=0,2,4$ we can observe stripes in the positive controls (shown by blue stars) above the diagonal for $x_i<x_d$ and below it for $x_i>x_d$, creating connections that will bring agents closer to the target. At $t=0$ the width of this stripe corresponds exactly to the confidence bound of the interaction function, as the control does not create edges that will not be used. Edges are removed to prevent individuals crossing the target, for example the red cluster visible at $t=4$ prevent agents' opinions overshooting. 

The following results provide an alternative characterisation of consensus at $x_d$ that demonstrates the two key challenges in this network control problem. Firstly we introduce some necessary notation. 

Denote the minimum and maximum opinions in the population by 
\begin{align*}
    x_m(t) = \min_{i\in\Lambda} x_i(t) \,,\quad
    x^M(t) = \max_{i\in\Lambda} x_i(t) \,. 
\end{align*}
For a given interaction function $\phi:[0,2]\rightarrow[0,1]$ we denote the set of roots of $\phi$ by 
\begin{equation}
    \mathcal{R}_\phi = \{ r \in [0,2] : \phi(r) = 0 \}.
\end{equation}
If $\mathcal{R}_\phi$ is empty then define $r^* = 2$, otherwise let $r^* = \inf (\mathcal{R}_\phi)$. 
    
\begin{proposition} \textup{(Proposition 2.1 in \cite{nugent2024steering})} \label{Propostion: x_d must remain in the opinion interval}
    Assume that the population reaches consensus. Then the population reaches consensus at a point $x_d\in[-1,1]$ iff $x_d\in [x_m(t),x^M(t)]$ for all $t\geq 0$. 
\end{proposition}

\begin{proposition} \textup{(Proposition 2.2 in \cite{nugent2024steering})}  \label{Proposition: gaps of size r^* cannot be closed}
    If $\phi$ is decreasing and there exists $i\in\Lambda$ such that $|x_{i+1}(0) - x_i(0)|>r^*$, then the population will not reach consensus. 
\end{proposition}

Together Proposition \ref{Propostion: x_d must remain in the opinion interval} and Proposition \ref{Proposition: gaps of size r^* cannot be closed} show that in order to achieve consensus at $x_d$ the population must be kept in one single opinion cluster, preventing any fragmentation, and the target opinion must remain within the range of current opinions. The first condition is particularly challenging as it is difficult to predict the clustering behaviour when considering non-linear opinion dynamics on a complex network, while the second condition show that a high level of precision is required as it is not possible to recover from overshooting the target opinion. This control problem is therefore challenging both theoretically, as many existing controllability results apply only to linear dynamics (see for example \cite{nozari2019heterogeneity,liu2018polarizability}), and numerically due to the rapidly growing dimension of the control space ($N^2$). For these reasons various extensions that have been considered for direct control (or similar problems in other interacting particle systems), such as sparse controls or partial information, remain open problems for network control.

\subsection{Discussion}

We complete this section with a discussion of the advantages and disadvantages of using ordinary and partial differential equations paired with control theory as a method to model opinion manipulation. 

As has been demonstrated by the various models in this section, a key advantage is the flexibility of this approach. If the effect of a control is known, such as altering network edge weights or directly impacting individuals' opinions, it is often straightforward to place a control in an appropriate location in an existing model to capture this effect. For example, if it were not possible to advertise to specific individuals but only to present a single opinion, then the individual $u_i$ term in \eqref{Eqn: directly affect dynamics} would be replaced with a single common control $u$. Moreover, the goals of a control can be clearly formulated through the cost functional. No prior knowledge or expectation of the behaviour of the control is required, only a decision regarding the placement of controls and structure of the cost functional. This allows optimal control to discover strategies that may not be obvious or immediately intuitive and thus improve our understanding of the system. 

However, such flexibility naturally also creates greater dependence on modelling choices. There may not be an obvious formulation or placement of a particular control. Consider for example the network control discussed in Section \ref{Section: Network Control}, in this setting it would also be possible to entirely replace the network $w_{ij}$ with a control $u_{ij}$, rather than placing the control in the derivative of $w_{ij}$. This leads to a new interpretation of the control, potentially new capabilities of the control and new considerations, such as whether this control should be restricted to be continuous in time or have some maximum rate of change. However, both replacing the network with a control and adding a control into existing network dynamics fulfil the criteria of including a control to model network manipulation. Just as the behaviour of the optimal control may not be obvious, the effect of these choices may not be obvious, and so the flexibility of this modelling approach is not purely advantageous. Modellers face a choice between presenting a more general framework and grounding problems in a specific context, but in either case should consider the robustness of control problems to their modelling choices, and ensure controls are added in a realistic way. 

Moreover, there are also choices in the construction of the cost functional. For example, \eqref{Eqn: directly affect cost functional} includes the difference from the target in the running cost (meaning it is placed inside the integral) rather than as a terminal cost (meaning it is evaluated only at time $T$). Including this as a running cost typically leads to faster convergence to the target, as there is no incentive for fast convergence when using only a terminal cost. However, this choice does not arise from any physical interpretation of the cost. It is physically meaningful to include the cost of controls (modelling for example the actual financial cost of advertising) but again a choice must be made regarding the exact form of this cost. A quadratic cost term is typically mathematically convenient, but this may not accurately reflect real costs in a specific scenario. Instead one could consider a cost functional of the form
\begin{equation}
    \mathcal{J}(x,u) = \int_0^T\left(\sum_{i=1}^N\frac{1}{2}(x_i - x_d)^2 + \frac{\nu}{2}\sum_{i=1}^N |u_i|\right)\, dt \,.
\end{equation}
This simple change from $u_i^2$ to $|u_i|$ can have a significant impact as it encourages the control to be sparse, meaning that a large number of controls may be set to $u_i=0$. This may more accurately reflect a scenario in which a controller wishes to influence as few individuals as possible. 

However, this change also significantly increases the complexity of the problem as the cost functional is no longer differentiable and so alternative optimisation techniques are required. The problem of sparse optimal control is addressed in the case of control through lying \cite{glendinning2025what} and in the case of the Cucker-Smale model for flocking, which is closely related to the HK formulation, in \cite{caponigro2015sparse,bailo2018optimal}. In \cite{fornasier2014mean}, Fornasier et. al provide a framework for more general cases of sparse optimal control problems with ODE constraints. Such problems are often already computationally intensive, with for example the network control problem requiring the repeated solving of $N^2$ coupled equations, which may limit their ability to provide insight as the problem must be solved for a range of initial conditions and parameters. 

Ultimately the goal of framing opinion manipulation as an optimal control problem is to understand the effectiveness and efficiency of a control mechanism. The cost of the optimal control provides a measure of the difficulty of a particular task. The control setup also allows one to consider the question of controllability or stabilisability \cite{sontag1998mathematical}. While a proof of controllability does not necessarily provide an explicit, efficient, or directly implementable control, it can be helpful to rigorously characterise the capability and limits of a particular control. For example, the direct control considered in Section \ref{Section: Directly Affecting Opinions} could be described as more powerful than the network control considered in Section \ref{Section: Network Control} as under this control any target opinion $x_d$ is globally asymptotically controllable, while this is not the case for network control. However, this assumes the controller can directly affect individuals' opinions, which may be less realistic than altering the balance of their social media content. 

While it may initially appear that modelling manipulation through controls assumes the presence of an external controller, the control through lying discussed in Section \ref{Sec: Liars} demonstrates that this is not the case. In particular, considering multiple interacting liars brings the control setup closer to that of game theory, whose techniques are discussed further in later sections.

\section{Binary Opinion Dynamics} \label{Section: Binary Opinion Dynamics}

In this section, we shift our attention to binary opinions; an alternative formulation of opinion dynamics, initially inspired by spin models in statistical mechanics \cite{liggett1985interacting}, typically with $x_i\in\{0,1\}$ or $x_i \in \{-1,1\}$, rather than the continuous opinions discussed in Section \ref{Section: Cts opinion dynamics}. In the absence of some (yet not all) of the analytical tools from continuous dynamics, problems are also studied from an algorithmic perspective. In particular, opinion manipulation is a prominent topic in the fields of algorithmic game theory, social choice, interaction protocols and mechanism design (see, e.g., \cite{CS404_guideBook} for an introduction). The emphasis shifts from whether a problem can be solved to how it can be solved efficiently, placing the focus on the design of the solution itself. From this angle, additional questions emerge, such as:
\begin{enumerate}
    \item What are the upper and lower complexity bounds of a given task, and are these bounds tight? \label{Q: Complexity}
    \item How robust is the algorithmic design? What guarantees does it offer regarding approximation error relative to the optimal solution and the rate of convergence (or termination time)? \label{Q: Convergence}
    \item What are the trade-offs between convergence speed and the size of agents' internal state spaces? \label{Q: Trade-off}
\end{enumerate}

In Section \ref{sec: intro binary opinions}, we present an overview of the main binary opinion-dynamics models we will discuss later on and introduce some relevant terminology. Specifically, in Section \ref{sec: influence max}, we consider controls that target agents directly, either by influencing their opinions or by modifying the links between them. Further, if one can manipulate specific agents, which ones should be chosen to maximise the effect? Section \ref{sec: coord and coop} and Section \ref{sec: exact majority} discuss how the influence maximisation techniques can be applied to two particular case studies. First, to nudge agents into cooperating, and second, to fulfil the task of inferring the majority opinion on the network from random local interactions.
Finally, we explore how manipulation can reveal hidden properties of the network in Section \ref{sec: inference}. 
We ask how hard it is to reconstruct the system's characteristics, with particular interest in the contrast between active and passive learning.

\subsection{Binary Opinion Dynamics Models}
\label{sec: intro binary opinions}

Binary opinion dynamics offer a minimal yet expressive framework for modelling agents' influence on one another. Let $L$ be the set of opinion labels, with $\vert L\vert=2$. Some standard choices, such as $L=\{0,1\}$ or $L=\{-1,1\}$, equip the system with arithmetic between labels, often interpreted as opposite or collapsing states. For example, having \emph{empty} versus \emph{full} spins \cite{liggett1985interacting}; holding \emph{for} or \emph{against} positions \cite{AAMAS24}; or being in a \emph{true} or \emph{false} state in circuits \cite{chistikov_convergence_2020}. Other models use entirely symbolic labels, such as \emph{red} and \emph{blue} ballots ($L=\{r,b\}$ as in \cite{grandi_complexity-theoretic_2025}) to allow for more nuanced interactions between agents known as population protocols \cite{ninjas2018}. The models we introduce in this section will have Markovian opinion updates. They only depend on the opinion states of the agents chosen to interact. Thus, when speaking of memory, we refer to the number of possible states an algorithm should account for to inform its output.

In the \emph{voter model}, each individual $i=1,\dots, N$ holds a binary opinion $x_i \in \{0,1\}$, and updates it by interacting randomly with others via a fixed network. In fields like automata theory or population protocols, the random interactions are thought of as being controlled by a \emph{scheduler} \cite{ninjas2018}. Networks can be represented either by an adjacency matrix $W$ or a graph $G=(V, E)$, composed of vertices $V$ and edges $E$ representing the agents and their connections, respectively. The influence weight between agents $i$ and $j$, $i\neq j$, is given by $w_{ij}$. To simulate the voter model in discrete time, at each time step, select one individual $i$ uniformly at random from the population, then choose a neighbour $j$ and adopt their opinion with probability
\begin{align*}
    p_{ij} = \frac{w_{ij}}{\sum_{k=1}^N w_{ik}} \,.
\end{align*}
Note that when all individuals hold the same opinion, the system is in a stable state. Weights may be binary, $w_{ij}
\in\{0,1\}$, to express the presence or absence of a link, or real-valued, $w_{ij}\in (0,1)$, representing the influence strength of a neighbour. Higher weights correspond to a higher likelihood of adopting that neighbour's opinion. 

The voter model can also be posed in continuous time by assigning a Poisson process $P_i$ to each individual and letting them update their opinion whenever an event occurs in $ P_i$. The resulting Markov process is studied extensively using duality and coupling in \cite{liggett1985interacting} for infinite populations. The papers discussed in this section, however, focus on finite populations. In either case, the system's behaviour depends heavily on the underlying network's topology.

\begin{figure}[ht!]
    \centering
    \includegraphics[width=0.75\linewidth]{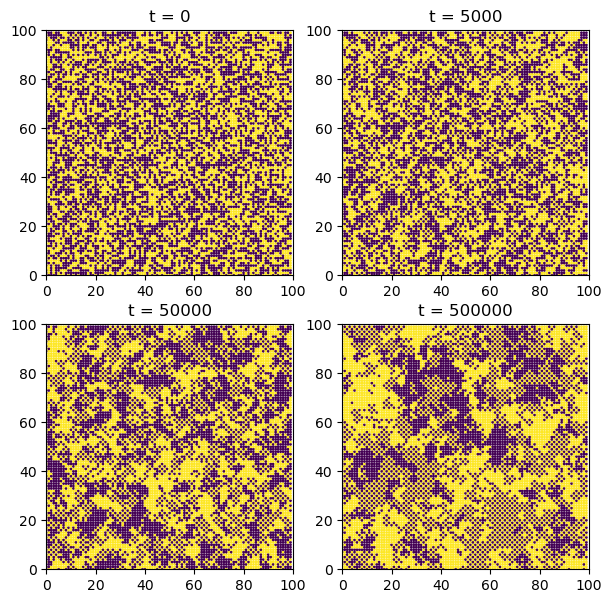}
    \caption{Realisation of the voter model with $N=1000$ individuals, showing the formation of spatial clusters. Agents are placed on a 2D lattice with connections to their four neighbours. Agents coloured yellow have opinion $0$ and agents coloured purple have opinion $1$.}
    \label{fig:VM}
\end{figure}

When working with a complete network, meaning $w_{ij}=1$ for all $i,j$, it is sufficient to track the number of individuals $n$ with opinion $0$ (or $1$). Since each interaction event only changes the opinion of one individual, $n$ makes jumps of size $1$. 
Denote by $r(n,n')$ the rate of jumps from $n$ to $n'$, then 
\begin{align*}
    r(n,n+1) &= \sum_{i,j=1}^N \mathbb{P}(x_i = 1) \mathbb{P}(x_j = 0) = \sum_{i,j=1}^N \frac{N-n}{N} \frac{n}{N} = (N-n)n \,, \\
    r(n,n-1) &= \sum_{i,j=1}^N \mathbb{P}(x_i = 0) \mathbb{P}(x_j = 1) = \sum_{i,j=1}^N \frac{n}{N} \frac{N-n}{N} = (N-n)n \,.
\end{align*}
Thus, $n\in\{0, N\}$ are fixed points corresponding to consensus states. When the network is complete (meaning $w_{ij}=1$ for all $i,j$), the probabilities $\mathbb{P}(x_i = 1) \mathbb{P}(x_j = 0)$ and $\mathbb{P}(x_i = 0) \mathbb{P}(x_j = 1)$ dependent on $i$ and $j$, so the rates above cannot be written in terms of $n$ only. A classical choice is to place individuals on a lattice with connections to their four (in 2D) nearest neighbours \cite{liggett1985interacting}. A realisation of the voter model with $N=1000$ individuals placed on a 2D lattice is shown at four time points in Figure \ref{fig:VM}. Agents coloured yellow have opinion $0$ and agents coloured purple have opinion $1$. As $t$ increases, the population forms into spatial patches, or clusters, in which individuals share the same opinion.

\emph{Majority Dynamics} defines a similar opinion-diffusion model, based on voting mechanisms. In the asynchronous case, the randomly selected individual $i$ updates their opinion to match the strict majority of their neighbours, with additional rules to break ties \cite{berenbrink_asynchronous_2022, Asynchronous_OD}. On directed graphs, only agents with incoming edges are considered for the vote; these agents are also known as the \emph{influencers} \cite{AAMAS24}.

In the synchronous case \cite{BredereckElkindManipulation, chistikov_convergence_2020, AAMAS24, EstradaChistikovPatersonTurrini2025}, all agents update their opinions simultaneously, meaning transitions occur between labelled social networks. Let $G=(V, E)$ be the graph representing a social network and let $\ell: V \to L$ be a \emph{labelling} function that assigns an opinion in $L$ to each agent in the network. Transitions are defined by the map $(G, \ell)\mapsto (G, \ell^+)$, where is given by:
    \begin{align}
        \ell^{+}(i) := \begin{cases}
           \ell(i)^{c} &\text{if } \; \vert\{j\in G_i : \ell(j)=\ell(i)\}\vert < \vert\{j\in G_i : \ell(j)\neq\ell(i)\}\vert,\\
            \ell(i) &\text{otherwise,}
        \end{cases}
        \label{eq: synchronous majority update}
    \end{align}
where $G_i$ is the set of influencers of agent $i$ and $\ell(i)^{c}$ is the opposite opinion of $\ell(i)$. 

Similar to the voter model, the behaviour of a system with majority dynamics depends heavily on the underlying network’s topology. For example, in directed acyclic graphs (DAGs), tree-like structures whose nodes act as either sources or sinks, the influence flows along well-defined paths. In fact, opinions are guaranteed to converge in at most the number of steps given by the length of the longest path \cite{chistikov_convergence_2020}. Further, any directed network can be reduced to a DAG by condensing its strongly connected components. To do this, each cycle is collapsed into a single ``supernode'' (see, e.g., Tarjan's algorithm \cite{Tarjan_alg}), and an edge is added between two supernodes whenever the original graph had an edge from one component to the other. The resulting condensation graph is a smaller DAG that preserves all reachability properties of the original network.

We refer to an external entity driving the manipulation as the \emph{campaigner}. Within a fixed opinion-diffusion system, the campaigner is endowed with some persuasive power; they may bribe agents to change opinions \cite{kempe2005influential}, add or remove links \cite{castiglioni_election_2020}, alter the update order \cite{BredereckElkindManipulation}, or use broadcast channels for advertisement \cite{Castiglioni2019ElectionMO, deligkas_being_2023}. Both the Voter Model and Majority Dynamics are examples of opinion-diffusion rules. They dictate how local opinions interact. This raises the question: what if the manipulation does not target the agents or the network, but the rules themselves?

\emph{Population protocols} offer a decentralised perspective on agent-based models \cite{ninjas2018}. They consist of $n$ identical finite-state machines (agents) governed by a deterministic transition function, activated by their random interactions moderated by the scheduler. The joint states of all agents constitute a \emph{configuration}, and a sequence of configurations evolving following the protocol's transition function forms an \emph{execution}. Moreover, an execution is said to be \emph{fair} when it does not avoid configurations that are reachable infinitely often. Equivalently, an execution is fair if every terminal state is reached with non-zero probability. That said, protocols are purposeful. They are sets of instructions designed to achieve a target task with high probability. The success of a protocol depends on whether executions converge to configurations that match an intended outcome, a property referred to as the protocol’s \emph{correctness}. Some prominent tasks in manipulating social network dynamics include those of \emph{cooperation and coordination} and \emph{exact majority}, which we will describe in further detail in sections \ref{sec: exact majority} and \ref{sec: coord and coop}, respectively. 

The exact majority task involves agents with binary states representing for and against opinions, and aims to have agents converge to a consensus on the initial global majority. Bear in mind that agents can only know the opinion of those whom they randomly bump into. Similarly, the cooperation and coordination task requires all agents to converge to the same opinion state. In this case, the goal is for agents to choose to cooperate, commonly in a \emph{prisoner's dilemma} setup, where, despite cooperation being the social welfare-maximising strategy, agents have selfish incentives to deviate. 

Population protocols exhibit a different form of manipulation; unlike the campaigner, they do not tamper with the agent's opinions, but instead change the local interaction rules to ensure that opinions eventually reach a desired state. This conforms to the branch of \emph{mechanism design} in distributed AI, the dual of game theory. It focuses not on how agents can optimally react to a given environment, but on how to design an environment so that, if agents act rationally, they exhibit the desired behaviour.  

\subsection{Influence Maximisation}
\label{sec: influence max}

The problem of how to maximise the influence that a control has over a population in order to cause the domination of a particular opinion has been the topic of much research in both continuous and discrete models. Initially, our control will come in the form of two stubborn agents who hold one of two binary opinions. Our problem is then reduced to considering which agents the controllers should be connected to (and what weight to give this connection) in order to induce a state that is favourable toward the controller's opinion. The controllers can be directly interpreted as social media bots, as explored in \cite{yadav2017influence}. Influence maximisation has been studied both from a view to characterise the complexity of the problem \cite{kempe_maximizing_2003, kempe2005influential} and to discover how algorithms can be used to explore strategies to maximise influence and uncover behaviours \cite{badawy2018analyzing, yadav2017influence}. In this section, we will consider a population with binary opinions $0$ or $1$, evolving according to the dynamics of the continuous voter model \cite{holley1975ergodic}, and seek to maximise the probability of each individual having a particular opinion. This problem is discussed in detail in \cite{brede2018resisting, romero2021shadowing, cai2025competitive, cai2022control}.

\begin{figure}[!ht]
\centering
\resizebox{0.4\textwidth}{!}{
\begin{circuitikz}
\tikzstyle{every node}=[font=\LARGE]
\draw (8.75,9) to[short, -o] (10.75,9) ;
\draw  (6.75,9) circle (2cm);
\node [font=\LARGE] at (10.75,8.5) {$1$};
\node at (7.25,10.25) [circ] {};
\node at (6.5,8.75) [circ] {};
\node at (6,9) [circ] {};
\node at (6.75,9.75) [circ] {};
\node at (7.75,8.75) [circ] {};
\node at (6.5,7.75) [circ] {};
\node at (5.5,8.75) [circ] {};
\draw (6,9) to[short] (6.75,9.75);
\draw (6.5,7.75) to[short] (5.5,8.75);
\draw (6.5,7.75) to[short] (6.5,8.75);
\draw (6.5,8.75) to[short] (7.75,8.75);
\draw (7.75,8.75) to[short] (6.75,9.75);
\draw (6.75,9.75) to[short] (7.25,10.25);
\draw (4.75,9) to[short, -o] (2.75,9) ;
\node [font=\LARGE] at (2.75,8.5) {$0$};
\end{circuitikz}
}
\caption{Description of the set-up for a system with $N$ agents and two controllers, one, labelled $0$, favours opinion $0$ and the other, labelled $1$, favours opinion $1$. The network of agents, shown in the centre of the circle, is connected with weights $w_{ij}$. Controllers $0$ and $1$ can be connected to any agent in the network, with weights $w_{0i}$ and $w_{1i}$, respectively.}
\label{Fig: influence maximisation set-up}
\end{figure}
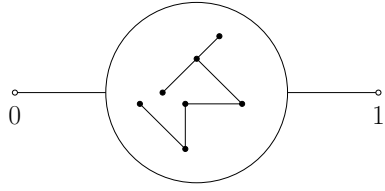

We consider the voter model outlined in Section \ref{sec: intro binary opinions}. In addition to our agents, suppose we have two external controllers corresponding to the two binary opinions $0$ and $1$. Controllers can freely decide which nodes they target with strengths $w_{0i}, w_{1i}\in\mathbb{R}^+$. For now, as discussed in \cite{romero2021shadowing}, we will assume only the controller favouring opinion $0$ is active, meaning that controller $0$ is seeking to determine the best way to distribute their link weights $w_{0i}$ in order to maximise the vote share for opinion $0$ amongst the population. We also place a budget constraint on the use of the control for opinion $0$, $\sum_iw_{0i}\leq\mathcal{B}_0$. In contrast, we treat link weights of the controller, favouring opinion $1$ $w_{1i}$, as being fixed and known by the active controller. We call the controller favouring opinion $1$ passive. Figure \ref{Fig: influence maximisation set-up} provides a diagram that helps explain this system.

We may study the opinion dynamics of this system by considering the evolution of the probabilities $p_i(t)\in[0,1]$ that a node $i$ holds opinion $0$ at time $t\geq 0$ \cite{mobilia2007role}. It is shown in \cite{romero2021shadowing} that the relevant rate equation in this case is
\begin{equation}
    \dot{p}_i = (1-p_i)\frac{\sum_jw_{ij}p_j + w_{0i}}{d_i + w_{0i}+w_{1i}} - p_i\frac{\sum_jw_{ij}(1-p_j)+w_{1i}}{d_i + w_{0i} + w_{1i}},
\end{equation}
where $d_i$ is the weighted degree of node $i$, i.e., $d_i=\sum_jw_{ij}$. The steady state of this system, $\mathbf{p}_*$, satisfies the equation $(\mathcal{L}+W_0+W_1)\mathbf{p}_*=\mathbf{w}_0$ \cite{masuda2015opinion,romero2021shadowing}. Here, $\mathcal{L}$ is the graph Laplacian of the underlying network, $W_0$ is a diagonal matrix with entries $w_{0i}$ and similarly for $W_1$. The vector $\mathbf{w}_0$ has entries $w_{0i}$. Hence, in order to maximise the probability of individuals holding opinion $0$, we choose our control weights $w_{0i}$ according to the following optimisation problem:
\begin{equation} \label{Eqn: Influence max utility function}
    \max_{W_0}\frac{1}{N}\mathbf{1}^T(\mathcal{L}+W_0+W_1)^{-1}W_0\mathbf{1},
\end{equation}
subject to $\sum_iw_{0i}\leq \mathcal{B}_0$ and $w_{0i}\geq 0$. This optimisation problem is concave, see \cite{romero2021shadowing}, hence convex optimisation techniques such as gradient descent may be applied.

Several strategies are observed for the active controller targeting opinion $0$. The first, referred to as \textit{shadowing}, is characterised by the active controller choosing to target the same nodes as the opponent controller. Analytic progress is possible by making an assumption that the node degree of ordinary agents is much larger than the weights attributed to the controllers, i.e., $(w_{0i}+w_{1i})\ll d_i$. Expanding about the resulting small parameter, authors in \cite{romero2021shadowing} determine optimality conditions of various orders. Shadowing is a first-order optimal response according to the resulting heterogeneous mean-field approximation to the opinion dynamics (see \cite{romero2021shadowing} for details). Surprisingly, given the nature of influence maximisation, the corresponding solution does not show any dependence on node degree $d_i = \sum_{j}w_{ij}$. In fact, the controller is most successful at maximising the first order cost when there is more deviation in the individual node allocations of the passive controller $w_{1i}$ from the mean $\frac{1}{N}\sum_jw_{1j}$. This means that when controller $1$ targets all ordinary nodes with the same amount of resource, the shadowing strategy of the active controller $0$ is less effective. The second strategy that will be discussed here is known as \textit{shielding}. A controller attempting to use a shielding strategy will target nodes that are direct neighbours of nodes being influenced by the passive controller. This is a second-order response to the optimisation problem with respect to the small parameter resulting from the assumption $(w_{0i}+w_{1i})\ll d_i$. It is noted that the effects of shielding and node degree are related in that nodes with high degree are more likely to be neighbours of nodes targeted by the passive controller.

We will now consider a variation of this model where both controllers $0$ and $1$ are active, hence influence maximisation is competitive. Clearly, this formulation is strongly linked to game theory approaches to considering opinion dynamics, discussed further in Sections \ref{sec: coord and coop} and \ref{sec: exact majority}. The controller favouring opinion $0$ seeks to maximise the vote share of option $0$ at some time $T>0$,
$$ S_0(T) = \frac{1}{N}\sum_{i=1}^Nx_i(T),$$
while the controller favouring opinion $1$ seeks to maximise $S_1(T) = 1-S_0(T)$. In \cite{cai2025competitive}, Cai et al. explore the problem of inter-temporal maximisation in a game-theoretical setting. We suppose we have given controls $w_{0i}(t)$ and $w_{1i}(t)$, such that
$$ w_{0i}(t) = \begin{cases}
    0 \quad & \text{for }0\leq t\leq t_0\\
    \tilde{w}_{0i} \quad &\text{for } t>t_0 
\end{cases}$$
where $\tilde{w}_{0i}$ is a given constant, typically $\tilde{w}_{0,i} = \frac{\mathcal{B}_0}{(T-t_0)N}$. The corresponding function $w_{1i}(t)$ is defined similarly. We are interested in identifying optimal start times for the control. Therefore, our minimisation problem is to find a $t^*_0$ and $t_1^*$ such that
\begin{subequations}\label{eqn: competetive influence max}
    \begin{align}
    t_0^* &= \arg\max_{0\leq t_0\leq T}\{S_0(T)\}\\
    t_1^* &= \arg\max_{0\leq t_1\leq T}\{1-S_0(T)\}
\end{align}
\end{subequations}
for $0\leq t_0\leq T$ where $T>0$ is an arbitrary time horizon. The existence of pure strategy Nash equilibria is theoretically guaranteed \cite{cai2025competitive}, however the computation through analytical means is demanding. Instead, the authors of \cite{cai2025competitive} utilise an iterative search algorithm, similar to that in \cite{bonomi2012computing}. The optimal starting times for the competitive controls can be explored numerically for different levels of heterogeneity of the underlying network. It is determined in \cite{cai2025competitive} that a controller with a budget advantage should activate their control earlier to ensure domination of vote shares. The model above can be easily adapted to give different nodes different start times for the allocation of control, i.e., a node $i$ is influenced by control $0$ at time $t_{0_i}^*$ for each $i=1, ..., N$. Indeed, it is also suggested that in regimes of sparse budget allocation that allocation should be strategic and the control should focus on low-degree nodes. 

Statistical inference methods for this type of model have also been explored \cite{cai2022control} as well as further influence maximisation questions such as how the control changes under the presence of biased voters (known as zealots) \cite{romero2020zealotry}. In a simpler case \cite{StubbornAgents} authors study the impact of stubborn agents in the voter model - agents who do not change their opinion in time but still influence others. Finally, questions of temporal influence maximisation have also been addressed in non-voter models, such as the SIS (susceptible and infected) model by Deligkas et al. \cite{deligkas_being_2023}, where the authors focus on the complexity and intractability of the resulting problem. 

Strongly related to influence maximisation, the problem of determining a minimum dominating set of a social network has been extensively studied \cite{sun2017dominating,nacher2016minimum,kempe_maximizing_2003}. For a graph $G=(V,E)$, the \textit{minimum dominating set} (MDS) of the graph is an optimised subset of nodes such that each node of the network either belongs to the subset or is adjacent to an element of the subset. Determining an MDS of a social network would be a cost-effective way to implement a manipulation strategy, since targeting only the members of the MDS would influence the entire population. While authors in \cite{nacher2016minimum} do not consider dynamics resulting from a certain selection of agents, the question of which nodes are the most efficient to target is still of interest for the study of manipulating opinion dynamics \cite{BredereckElkindManipulation, AAMAS24}. 

In \cite{sun2017dominating}, the authors discuss an alternative definition of a dominating set that accounts for the graph's structure. A \textit{structure-driven minimum dominating set} (SD-MDS) assumes that a driver node (a node belonging to the SD-MDS) can control a non-driver node if the link between them belongs to a cyclic structure. Considering the minimum dominating set is highly applicable to social networks, where a cyclic structure can be interpreted as a friendship group, where non-connected nodes still influence one another through friend-of-a-friend interactions. If $G=(V, E)$ is a graph with suitable structure, controlling agents within an SD-MDS is a much more cost-effective way to control the entire population.

Another model tackling problems similar to influence maximisation is discussed in \cite{kozitsin2022general,kozitsin2022optimal,kozitsin2024optimal}, where agents may have opinion values taken from a finite set $L=\{a_1, ..., a_m\}$, as seen in the labelling model outlined in Section \ref{sec: intro binary opinions}. We introduce \textit{stubborn agents} (who do not change their opinions in time) to control the opinions of the population. Similarly to in Section \ref{Section: Cts opinion dynamics}, we denote the opinion of agent $k$ at time $t\geq0$ by $x_k(t)$. If there are $M$ stubborn agents, we can split this population into stubborn agents targeting ordinary agents with opinion $a_i$ for each $i=1, ..., m$ to demonstrate personalised influence strategies, often used by bots. In the mean-field limit, detailed in \cite{kozitsin2024optimal} the dynamics of the proportion of normal agents with opinion $a_i$, denoted $y_i$ for $i=1, ..., m$ is given by
\begin{equation}\label{Eq: dynamics controlled influence maximisation}
    \frac{dy_i}{dt} = \sum_{r=1}^m\sum_{s=1}^my_r(t)\left[y_s(t)p_{r,s,i} + u_s(t)p^I_{r,s,i}\right] - y_i(t).
\end{equation}
Here, 
\begin{align*}
    p_{r,s,i} &= \mathbb{P}(x_k(t+1) = a_i|x_k(t) = a_r,x_l(t) = a_s), \\
    p_{r,s,i}^I &= \mathbb{P}(x_k(t+1) = a_i|x_k(t) = a_r,v_l(t) = a_s),
\end{align*}
$u_l(t)$ is the proportion of stubborn agents with opinion $a_l$ at time $t\geq 0$ and $v_l(t)$ is the opinion of stubborn agent $l$.

We suppose that a concerned person is interested in influencing the population on some time interval $[0, T]$ through adjusting the opinions of stubborn agents. Similarly to Section \ref{Section: Cts opinion dynamics} we may then formulate the following optimal control problem by seeking to minimise
\begin{equation*}
    \mathcal{J}(y,u) = K\int_0^T\sum_{i=1}^ms_iy_i(t)\, dt + \sum_{i=1}^ms_iy_i(T)
\end{equation*}
where here $s = [s_1, ..., s_m]^T$ is a non-negative vector of weights, assumed to be given, informing the desired opinion distribution. Furthermore, $K>0$ encodes the relative importance of the final opinion distribution. We solve our optimal control problem with respect to the dynamics given in equation \eqref{Eq: dynamics controlled influence maximisation} and the constraints $\sum_{i=1}^mu_i(t) = \frac{M}{N+M}$ and $u_i\geq0$ for $i=1, ..., m$. Analytic results on existence of optimal controls and numeric methods for determining controls are explored in \cite{kozitsin2024optimal}. It is conjectured based on numerical results that the control should focus its efforts on large opinion clusters of a particular opinion type, similar to a bang-bang control. In \cite{kozitsin2024optimal}, a more general framework is considered where normal agents are further divided into being of multiple type. Different sparse network structures such as Barabasi-Albert and Watts-Strogatz graphs are also considered.

Influence maximisation is a highly relevant and mathematically interesting problem inspired by modelling manipulation on a social network. The problem of when to trigger a control \eqref{eqn: competetive influence max}, outlined in \cite{cai2025competitive}, is important for understanding the behaviour of bots on social media for example. This is an area that has not been well explored in continuous opinion dynamics thus far. Competitive influence maximisation studies a similar framework to game theory which will be the focus of the following few sections.

\subsection{Cooperation and Coordination}
\label{sec: coord and coop}

We consider the cooperation and coordination problem in this section. The goal is to manipulate intrinsically self-interested agents into adopting the social-welfare-maximising strategy. Perhaps the most famous setup for this paradox is the prisoner’s dilemma, in which two agents can either incriminate the other to secure a better sentence or cooperate and remain quiet so neither goes to prison. The caveat is that if both defect by snitching on one another, they face the worst outcome together. A classical stage of the \emph{tragedy of the commons}. In real-world scenarios, cooperation would look like contributing to a shared resource fund to avoid a collective catastrophe (e.g. investments in green energy \cite{smirnov_collective_2019}) or incurring a cost in exchange for social services (e.g. taxes \cite{gottlieb_tax_1985}). At the same time, defection conforms to a free-rider effect \cite{free_rider}.

One possible solution, evident in real life, is to allow agents' connections to react to their behaviour. In modelling terms, this corresponds to a dynamic network. We introduced this in the context of continuous opinion dynamics in Section \ref{Section: Network Control}, but it has also been studied in the context of social dilemmas. The stream of papers \cite{grandi_complexity-theoretic_2025, BredereckElkindManipulation, auletta_complexity_2019, castiglioni_election_2020, li_evolution_2020} examines the computational aspects of exploiting the majority opinion dynamics, as in \eqref{eq: synchronous majority update}, using interventions such as opinion transformation \cite{kempe_maximizing_2003} and edge removal or addition \cite{harrell_strength_2018}, but now targeted at the specific task of fostering cooperation. Cooperation seems to emerge when agents trust that they are matched with other cooperative agents and when the rate at which they can sever links is sufficiently fast \cite{li_evolution_2020, harrell_strength_2018}. Agents would rarely break links with cooperators. Hence, a natural phenomenon in both real-life and simulated experiments \cite{BaraTurriniAndrighetto23, bravo_trust_2012, rand_dynamic_2011} is that dynamic networks evolve toward a structure with cooperators at the core and defectors in the periphery. Typically, cooperators are identified by high reputation or high degree, which, in turn, drives preferential attachment by other agents \cite{Preferential_Attachment, santos_picky_2020}. In continuous opinion dynamics, influential nodes are identified using time-varying masses \cite{piccoli2021control}.

\begin{figure}[ht]
\centering
\includegraphics[width=0.7\textwidth]{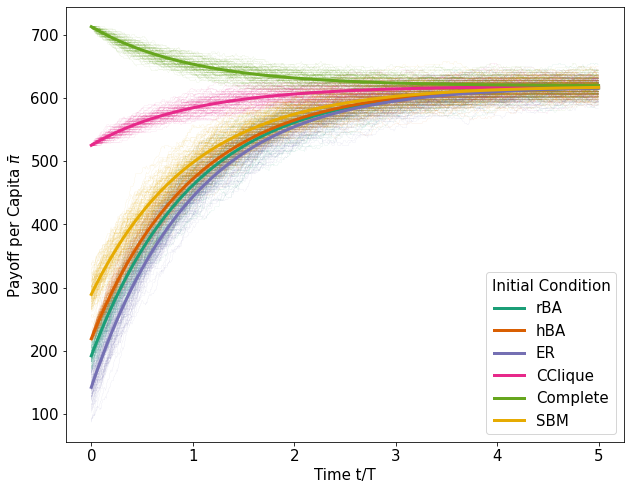}
\caption[Enabling imitation-based cooperation]{Evolution of payoff per capita as agents play a Prisoner’s Dilemma with fixed strategies and cooperators have extreme popularity. Faint lines indicate individual runs, while the bold lines represent the CANDY results. Time has been rescaled by the number of dyads, $T=N(N-1)/2=190$, so that every $T$ time-steps on average all pairs have been updated once. All payoffs per capita converge regardless of the initial condition (i.e., the network generation model) as all networks develop core-periphery structures. Reproduced with permission from \cite{BaraTurriniAndrighetto23}.}
\label{fig: Bara Imitation}
\end{figure}

Unlike Section \ref{Section: Network Control}, where the aim was to find the optimal control of an existing system, the goal here is to construct the network dynamics themselves in such a way that encourages cooperation. To design the payoff that triggers agents to become cooperation hubs, authors in \cite{BaraTurriniAndrighetto23} propose a graph-theoretic model of extreme popularity, in which only cooperators are befriended, and only defectors are unfriended. The respective rates of agent $i$ generating, or breaking, a link with agent $j$ are
\begin{align}
    g_{ij}(0,1)= \frac{1}{\tau_g}\frac{s_i + s_j}{2} \quad \text{and} \quad g_{ij}(1,0)= \frac{1}{\tau_g}\left(1 - \frac{s_i + s_j}{2}\right),
    \label{eq: make or break ties}
\end{align}
where $s_i, \; s_j \in \{0, 1\}$ are the strategies of agent $i$ and $j$, with $0$ to defect and $1$ to cooperate, and $\tau_g$ representing any graph-theoretic timescale. In particular, $\tau_g$ represents how the opinion updates sync up with the link changes of the agents. Figure \ref{fig: Bara Imitation} shows this protocol in action, with agents evolving towards the same strategy regardless of the initial network structure. 

\begin{remark}
    The emergence of cooperation is not exclusive to discrete opinions. In the continuous case, \cite{nugent2024steering} showed that the optimal update frequency is key to the emergence of cooperation. In fact, the optimal rewiring frequency depends on the rate of change of the underlying opinion dynamics. Classical notions such as stability, robustness, and timescale separation persist even when interaction protocols are affected by a feedback loop, that is, a learning-based environment \cite{recht_tour_2019}.
\end{remark}

Building on this foundation, \cite{chen_general_2024} proposes a general control-theoretic framework for reinforcement learning that establishes convergence and optimality guarantees through analogues of the Bellman operator used in Q-learning (see, e.g., \cite{CS404_guideBook}). Cooperation and coordination emerge as the learning dynamics that implicitly define the feedback laws shaping agent interactions over time. In parallel, \cite{ramadan_actively_2024} frames reinforcement learning as a stochastic optimal control problem in which agents actively trade off caution (\emph{exploitation}) and uncertainty (\emph{exploration}). In dynamic social networks, this manifests as agents aiming to strike a balance when adjusting their connectivity and trust in their neighbours in response to unexpected outcomes. A strategy proven to experimentally drive cooperation is the \emph{Out-for-Tat}, whereby agents break ties with each other once either agent defects \cite{PaoloOutForTat}. Together, these works suggest that fostering cooperation can be understood as a control problem in which both strategies and interaction protocols are co-designed with learning dynamics to regulate uncertainty.

Ultimately, cooperation and coordination arise not from static incentives or fixed network structures, but from an interplay between agent behaviour, learning dynamics, and adaptive connectivity. Dynamic networks provide a natural mechanism for reinforcing or weakening trust, reputation, and influence over time, while control-theoretic and reinforcement learning approaches offer formal tools for analysing and designing these processes. This perspective motivates studying endogenous network dynamics as part of the control loop itself, rather than as an external entity that sets the stage for principled methods to induce socially desirable outcomes in populations of self-interested agents.

\subsection{Exact Majority}
\label{sec: exact majority}
We consider the exact majority problem in this section. Agents hold an initial opinion and aim to infer the global majority through random interactions with other agents. Therefore, we are in the realm of population protocols mentioned in Section \ref{sec: intro binary opinions}. Our goal in this section is to show that by manipulating interaction rules a protocol can (or cannot) yield formal guarantees, using the exact majority task as a case study.

Formally, a population protocol is a tuple $\mathcal{P}=(Q, T, I, O)$, where $Q$ is a finite set of states, $T\subseteq Q^2\times Q^2$ encodes the transition rules, $I$ contains possible initial states and $O$ is a function $O: Q \to \{0,1\}$ that yields the binary output of the system. Transitions have the form $((p,q), (p',q'))\in Q^2\times Q^2$ to indicate that if two agents interact, one in state $p$ and another in state $q$, then they would change to states $p'$ and $q'$, respectively. A configuration $C\in \mathbb{N}^Q$ assigns a number of agents to each admissible state. If $C$ remains unchanged by the protocol, then it is \emph{terminal}. Otherwise, the possible subsequent configurations are said to be \emph{reachable} and denoted with $C \overset{*}{\to} D$, for some $D\in \mathbb{N}^Q$. An execution induces an infinite configuration's sequence $C_0 \overset{t_1 \;}{\to} C_1  \overset{t_2 \;}{\to} C_2 \dots $ for transitions $t_1, t_2, \dots \in T$, and consequently, it is fair if $\forall D$, if $\{i\in \mathbb{N}: C_i \overset{*}{\to} D$\} is infinite, then $\{i\in \mathbb{N}: C_i=D\}$ is also infinite. An \emph{outcome} of a configuration is determine by consensus, so that $O(C)=b \in \{0,1\}$ if $O(q)=b$ for every $q\in Q$ and $\bot$ otherwise. A configuration is stable if $O(C)\in \{0,1\}$ and if $\forall D,$ if $C \overset{*}{\to} D$, then $O(D)=O(C)$. That is, the configuration reaches consensus and stays there. A protocol $\mathcal{P}$ \emph{computes} a predicate $\varphi: \mathbb{N}^I\to \{0,1\}$ if $\forall C\in I$, every fair execution from $C$ leads to a stable configuration $D$ such that $O(D)=\varphi(C)$. Therefore, a protocol is correct with respect to a predicate if the predicate can predict the outcome of stable configurations. Naturally, a first research question is whether a protocol is correct under fair executions for a given task, and then, ask whether the protocol can be executed efficiently.

Authors in \cite{ninjas2018} draw an analogy between these interacting protocols' properties and a group of ninjas running in the dark who must decide whether to attack a fortress. Each ninja knows that others exist, while oblivious of how many ninjas there are in total or whether they support or are against the attack. Each ninja holds their opinion about whether they want to attack the fortress or not. Yet, they are instructed to follow the majority.

The problem, then, is that ninjas need to make their decision on the global majority whilst being limited to the encounters with the ninjas they randomly bump into. Thus, how informative such local interactions can be \cite{cord-landwehr_network_2015}? In real-world scenarios, this resembles the echo chamber effect on social media \cite{bara_predicting_2022} leading to a majority illusion among agents \cite{grandi_complexity-theoretic_2025}.

A protocol dictates how information is exchanged during collisions, and the random encounters correspond to the scheduler selecting two agents to interact. Authors in \cite{ninjas2018} extend their ninjas analogy to formalise notions of fairness, correctness, and convergence by mechanically verifying the majority protocol. They demonstrate that altering the local interaction rules can determine whether the protocol is correct and accelerate the convergence to consensus. After all, ninjas must decide before dawn. 

We summarise the protocols from \cite{ninjas2018} in Table \ref{tab: Three protocols}. Agents can be either \emph{Active} (A) or \emph{Passive} (P), and hold a belief regarding the majority’s stance towards an attack: \emph{Yes} (\textbf{Y}), \emph{No} (\textbf{N})\footnote{Do not confuse with population size $N$.}, or (only in Protocol 3) a \emph{Tie} (\textbf{T}). Generally, active agents impose their opinions on passive agents. Protocol 1 is incorrect when a tie occurs, as agents never reach consensus. Protocol 2 is correct, with agents reaching a consensus not to attack in the event of a tie. However, when the difference between \emph{Yes} and \emph{No} votes is minimal, convergence is slow: agents in the majority must first convert the minority into a passive–tied state, so that the last few active agents can transform them one by one into the true majority. Protocol 3 avoids these caveats, albeit at the cost of a larger state space, as agents must consider a tied voting outcome.

\begin{table}[ht]
\caption{Sets of rules from \cite{ninjas2018} for the task of computing the Exact Majority.}\label{tab: Three protocols}
\begin{tabular}{c|c|c}
\textbf{Protocol 1\footnotemark[1]} & \textbf{Protocol 2\footnotemark[1]}  & \textbf{Protocol 3\footnotemark[2]} \\ 
\midrule
$(A_\textbf{Y}, A_\textbf{N}) \mapsto (P_\textbf{N}, P_\textbf{N})$  & $(A_\textbf{Y}, A_\textbf{N}) \mapsto (P_\textbf{N}, P_\textbf{N})$  & $(A_\textbf{Y},A_\textbf{T}) \mapsto (A_\textbf{Y},P_\textbf{Y})$\\
$(A_{\alpha}, P_{\beta}) \mapsto (A_{\alpha}, P_{\alpha})$ & $(A_{\alpha}, P_{\beta}) \mapsto (A_{\alpha}, P_{\alpha})$ & $(A_\textbf{Y},A_\textbf{N}) \mapsto (A_\textbf{T},P_\textbf{T})$\\
& $(P_\textbf{N}, P_\textbf{Y}) \mapsto (P_\textbf{N},P_\textbf{N})$ & $(A_\textbf{T}, A_\textbf{N})\mapsto (A_\textbf{N}, P_\textbf{N})$ \\
&  & $(A_\textbf{T}, A_\textbf{T})\mapsto (A_\textbf{T},P_\textbf{T})$ \\
&  & $(A_{\alpha}, P_{\beta}) \mapsto (A_{\alpha}, P_{\alpha})$
\end{tabular}
\footnotetext[1]{ For $\alpha,\beta\in \{\textbf{Y},\textbf{N}\}$.} 
\footnotetext[2]{ For $\alpha, \beta \in \{\textbf{Y},\textbf{T},\textbf{N}\}$.}
\end{table}
Correctness alone does not guarantee implementability; an algorithm may yield the correct outcome yet require excessive time or space. Similar to the cost functionals used in optimal control, computational cost and memory convey a similar preference for efficient manipulation. When speaking about efficiency, one may either present an algorithm that guarantees the task can be completed with high probability within a prescribed time and space, or, conversely, prove that no algorithm can achieve this with fewer resources without breaking widely held hardness assumptions. 

In \cite{2017LeadersPopulationProtocols}, the authors propose the Split-Join algorithm for the Exact Majority problem, which uses $\mathcal{O}(\log^2 N)$ states per node. Each node stores an integer value, with the sign indicating whether it supports the majority. At each interaction, nodes average and normalise their values. So the total sum across all nodes matches the sign of the true majority, thus preserving correctness. The algorithm is guaranteed to converge to the correct majority within $\mathcal{O}(\log^3 N)$ parallel time using $\mathcal{O}(\log^2 N)$ states. Whether such a rolling-sum estimate is the most effective way to track the exact majority once there is online or incomplete information remains an open question. Similar tasks also aim to extend the Exact Majority to consider pluralism \cite{angluin_computational_2007, angluin_computation_2006, angluin_stably_2006, angluin_fast_2008}. 

The Exact Majority proves how additional rules can speed up convergence. Further, agents with an alternative set of rules can also reduce the space required. For example, in a more general task, \cite{angluin_fast_2008} showed that computing Presburger-definable predicates (i.e., predicates constructed using only linear inequalities and modular arithmetic) can be achieved an expected convergence time of $\mathcal{O}(N^2 \log N)$ when there's a homogenous set of rules for all agents, but that can then be improved to $\mathcal{O}(N\log^5 N)$ in presence of leaders (heterogeneous sets of rules) \cite{angluin_computation_2006, 2017LeadersPopulationProtocols}.

\subsection{Network Inference}
\label{sec: inference}
So far, most work on opinion manipulation has assumed a known network structure when identifying the most influential agents or links \cite{BredereckElkindManipulation, kempe2005influential, chen_efficient_2009, leskovec_cost-effective_2007, tang_influence_2015, kempe_maximizing_2003, castiglioni_election_2020}. These approaches address the questions of what, when, and who to influence for maximum impact. Yet, they overlook what a network's response to manipulation reveals about its structure. In this section, we draw attention to the converse problem and to how manipulation can be used to infer a network's hidden properties \cite{gomez-rodriguez_inferring_2010, AAMAS24, qiu_learning_2024, EstradaChistikovPatersonTurrini2025}.

Take the two underlying opinion dynamics described in \cite{kempe_maximizing_2003}: the Independent Cascade and the Linear Threshold. These are dual representations of the same underlying diffusion process that defines progressive contagions on 
agents with binary states. Namely, active and inactive (as initially proposed), susceptible and infected (in epidemiology contexts \cite{leskovec_cost-effective_2007}), or for and against (when speaking of opinions \cite{AAMAS24}). Note that the states are not necessarily symmetrical: once an agent becomes active, it cannot be deactivated; infected agents remain infected; and once an agent changes their opinion, there is no turning back. In the Independent Cascade, whenever an agent adopts a new opinion, it gets a chance to influence each of its neighbours. Upon success, the neighbour also adopts this new opinion; on failure, the edge is removed, and no further attempts are made. 

Formally, let $G=(V, E)$ be any directed graph with at most $k>0$ edges (i.e., $\vert E \vert \leq k$) and $N=|V|$. A cascade $c$ emerges as the opinion spreads from neighbour to neighbour across the social network and is stored as a hitting time vector $c:=[t_1, \dots, t_{n}]$ containing the time each agent adopted the target opinion, with the convention that $t_i=\infty$ if $i\in V$ is not reached. Denote by $\mathcal{T}_c(G)$ the set of all possible propagation trees in $G$ for a cascade $c$. Multiple cascades come together into a contagion trace $C$. Further, let $F_C(G)$ be the log-likelihood of $G$ admitting all the cascades in $C$. Thus, inferring the most likely network boils down to finding $G^*$ that induces $\mathcal{T}_C(G^*)$ with the most trees. We obtain the optimisation problem

\begin{align}
    G^* = \argmax_{G=(V,E) : \vert E \vert \leq k} F_C(G)= \sum_{c\in C} \max_{T\in \mathcal{T}_{c}(G)} \sum_{(i,j)\in E_T} w_{ij},
    \label{eq: network inference optimisation}
\end{align} 
where $w_{ij}$ is a nonnegative weight representing the improvement in log-likelihood of edge $(i, j)\in E_T$ under the most likely propagation tree $T=(V_T, E_T)$ contained in $G$. Without any further manipulation of the system, authors in \cite{gomez-rodriguez_inferring_2010} proved that the problem in \eqref{eq: network inference optimisation} is NP-hard due to the combinatorial explosion in the possible edge configurations. As an alternative, they proposed an approximation algorithm that achieves at least $63\%$ of the optimal value. 

For the \emph{linear threshold}, \cite{AAMAS24} considers a campaigner who aims to learn the structure of a social network by influencing agents' opinions and observing the subsequent opinion updates. In their setup, agents adopt the opinion of the strict majority of their influencers as described in \eqref{eq: synchronous majority update}. Further, the campaigner keeps track of an observation and intervention budget. The first refers to the number of opinion updates they can observe, and the second to the number of opinions the campaigner can directly change. In contrast to \cite{gomez-rodriguez_inferring_2010}, this approach yielded a positive result, proving that the network influence problem becomes tractable once manipulation is allowed. The authors provided polynomial upper bounds on the necessary budgets to guarantee exact learning of the underlying network. To do so, they proposed an algorithm for learning networks in general, which was subsequently refined for cliques (see Table \ref{tab: Summary of upper bounds}). 

\begin{figure}[ht]
\centering
\includegraphics[width=0.8\textwidth]{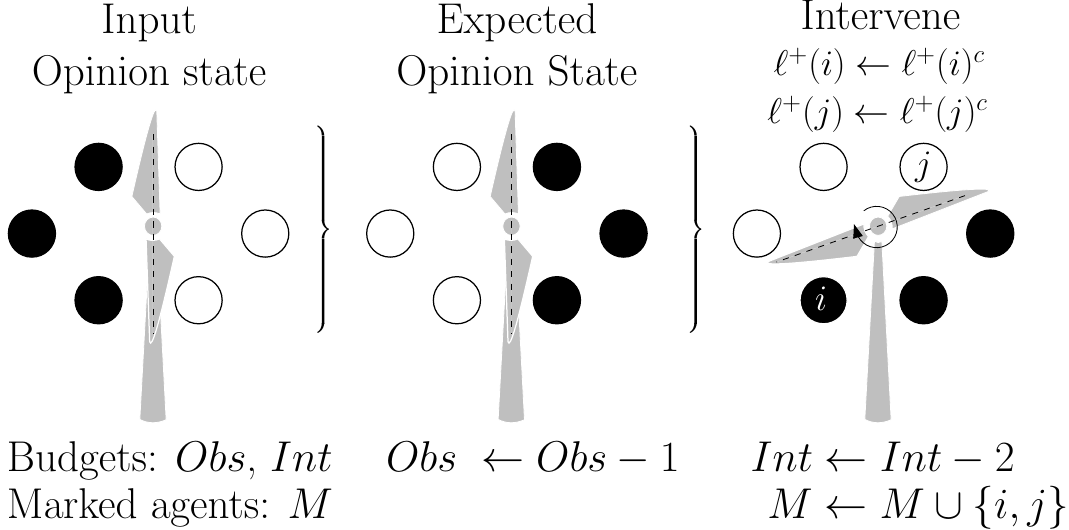}
\caption[Graphical intuition of \emph{The Windmill}]{Graphical intuition of the \emph{Windmill Algorithm}, reproduced with permission from \cite{AAMAS24}. Agents' opinions are depicted with \opinion{black} and \opinion{white}. The campaigner’s observation and intervention budgets, $Obs$ and $Int$, decrease as they are used. The set of marked agents $M$ keeps track of which agents have been intervened on.}
\label{fig: Windmill algorithm}
\end{figure}

\begin{table}[ht]
\caption{\centering Upper bounds on the campaigner's budget for different network-inference tasks under Majority Dynamics.}\label{tab: Summary of upper bounds}
\begin{tabular}{c|c|c}
    \textbf{Learning task} & \textbf{\renewcommand{\arraystretch}{1}\begin{tabular}[c]{@{}c@{}}Observation \\ budget \end{tabular}} & \textbf{\renewcommand{\arraystretch}{1} \begin{tabular}[c]{@{}c@{}}Intervention \\ budget \end{tabular}} \\
    \midrule
    Learn any network $G$& $\oO(N^2)$ & $\oO(N^3)$ \\
    Identify an odd clique & $\oO(N)$ & $\oO(N^2)$ \\
    Identify an even clique & $\oO(N)$ & $\oO(N)$
\end{tabular}
\end{table}
The best-case scenario occurs on cliques formed by an even number of agents, which can be inferred using linear resources. The campaigner is familiar with the opinion-diffusion protocol and uses symmetry-breaking arguments to improve the inference bounds through what they call the \emph{Windmill Algorithm}. Intuitively, when an even clique faces a 50/50 opinion split, all agents change their opinions. The campaigner spends their intervention budget to prompt these updates in search of a counterexample to the belief that the network is a clique. Yet, they cannot be sure the network is a clique from a single query. So, they need to conveniently arrange the agents in a circle, and intervene in pairs with opposite opinions, rotating like a windmill (as illustrated in Figure \ref{fig: Windmill algorithm}). Only cliques can complete a full rotation without running into any inconsistencies \cite{AAMAS24}.

Manipulation makes inference tractable by allowing the campaigner to choose which opinion updates to observe. This is often referred to as \emph{active learning} in distributive AI, where the learning algorithm selects a specific instance to query. In fact, once the campaigner loses their intervention power and opinions are drawn at random, constituting a \emph{passive learning} framework, inference becomes hard. Indeed, recent studies \cite{EstradaChistikovPatersonTurrini2025} and \cite{qiu_learning_2024} show that even attempting to learn a network \emph{probably approximately correctly} (PAC), see e.g., \cite{kearns_introduction_1994}, is NP-hard whenever the opinion dynamics obey any threshold-based rule. 

\subsection{Discussion}

We conclude this section with a discussion of the advantages and disadvantages of using algorithms to model binary opinions in the context of opinion manipulation. An advantage is the clarity and simplicity afforded by restricting the opinion space to two states. Several real-world scenarios naturally align with a binary formulation, allowing these models to capture essential system features without requiring additional memory. The restricted state space often leads to cleaner probabilistic descriptions. With only finitely many possible outcomes, it allows for easier characterisation of transition probabilities, convergence estimates, and analysis of random processes that would become unwieldy in a continuous setting.

Still, it is possible to borrow from the continuous toolkit. Transitions in the voter model mirror the adjacency matrix of an underlying Markov process, which allows the use of classical analytical techniques. Although the present work does not explicitly rely on duality theory, the correspondence between network topology and transition behaviour is a recurring theme in the literature and often shapes the available analytical approaches. Conversely, the emergence of specific network structures in dynamic networks, such as the core–periphery patterns observed for cooperation and coordination in Section \ref{sec: coord and coop}, is not unique to binary opinions, as it also appears in continuous settings \cite{nugent2024steering, nugent2023evolving}.

However, as with any framework, working with binary opinions entails modelling choices. It inevitably restricts the model’s expressiveness in settings involving graded beliefs, multidimensional preferences, or evolving internal states. Likewise, the formulation of control mechanisms changes substantially when only two states are available. While the cost of manipulation often admits a simpler form, this does not imply that it is unique or ``correct.'' For example, a control acting directly on opinions represents a different mechanism from one that modifies the underlying interaction pattern; or optimising controls within a fixed interaction structure versus optimising the way a network is rewired, even though they may be admissible modelling choices for similar tasks. There is a tension between generality and specificity. All these perspectives satisfy the formal requirement of adding a control, yet they lead to different capabilities and, crucially, different feasibility assumptions.

Beyond these structural limitations, additional considerations arise when formulating objectives. Approximate algorithms allow one to forgo complete accuracy in exchange for improved time efficiency. Choices regarding fairness assumptions, convergence guarantees, or the structural properties of the underlying interaction graph can affect both the tractability and the interpretation of the model. Algorithm design is more open-ended than an optimal-control approach. When a problem is computationally HARD\footnote{ HARD as a complexity class.}, algorithms do not necessarily yield definitive answers, but rather benchmarks. For instance, \cite{kempe_maximizing_2003} laid the groundwork for the influence maximisation problem discussed in Section \ref{sec: influence max}, and subsequent work \cite{leskovec_cost-effective_2007, chen_efficient_2009, tang_influence_2015} progressively developed algorithms that advanced the state of the art.

Differences may also arise when understanding the manipulation problem from a Game Theory or Mechanism Design perspective. The first line of work focuses on best-response behaviour and the rate of convergence to a Nash equilibrium. Instead, the latter considers whether it is possible to construct a set of incentives that ensures agents behave in a prescribed manner. Population protocol models sit at the intersection of these approaches. They aim to design rule sets that drive the system towards fulfilling specific tasks (e.g., exact majority or cooperation and coordination). Yet, lower bounds show that even carefully designed rules are subject to fundamental spatio-temporal constraints. In particular, when the number of internal states is severely limited, the expected convergence time cannot improve beyond $\Omega(n/\text{polylog} (n))$ \cite{ninjas2018}. This highlights a gap between what can be achieved through local manipulation and what is theoretically optimal. 

Viewing manipulation as a form of system design raises interpretative questions about the nature of control. Some models posit a central authority that directly shapes interactions or allocates incentives, while others assign the same mathematical structure to strategic behaviour emerging from decentralised agents. Although this broadens the modelling scope, it also complicates analysis, as system behaviour may become sensitive to assumptions about rationality, information availability, and timing.

At the same time, this perspective affords considerable flexibility. Once the behaviour of a baseline system is understood, additional layers of control can be introduced to enhance robustness or induce desired outcomes. This layered approach is illustrated in Section \ref{sec: exact majority}, where augmenting population protocols with additional rules improves their robustness. In this sense, a manipulated system can itself become the setting for further manipulation, enabling iterative refinement and resilience in the face of increasingly complex dynamics. 

Ultimately, binary opinion models offer an analytically convenient approach to studying influence and control, but their simplicity comes at a cost. As with any flexible modelling framework, the modeller must balance generality with interpretability, remaining aware of the consequences of their design choices and considering whether the assumptions underlying a particular formulation are suited to the phenomenon of interest.

\section{Conclusion}\label{Section: conclusion}

In Section \ref{Section: Cts opinion dynamics} and Section \ref{Section: Binary Opinion Dynamics} we have introduced the key approaches in modelling opinion formation in continuous and binary opinion spaces respectively, demonstrated methods to incorporate manipulation, and discussed their various advantages and disadvantages. We conclude here by highlighting the connections between some of these approaches, before taking a moment to give the authors' own perspectives on the impact and importance of this field. 

All the methods of modelling opinion manipulation discussed in Section \ref{Section: Cts opinion dynamics} and Section \ref{Section: Binary Opinion Dynamics} capture in some way the objective of the manipulation. This can typically be separated into two parts, the possibility of achieving some target and a measure of efficiency. Here we see parallels between two rather different approaches to modelling manipulation: the addition of a control to an existing model and the construction of an interaction protocol. In both settings one must begin by establishing the limits of manipulation. In control theory this is done by deciding at what location in an existing model a control will be introduced. In protocol design this instead appears through the requirement for fairness. We may then consider if a target can be achieved within these limits: in control theory this comes through a proof of controllability (alternatively observability or stabilisability, depending on the context) while in protocol design the goal is to show the correctness of a proposed protocol. Finally efficiency in control theory is measured through a cost functional (such as \eqref{Eqn: directly affect cost functional}, \eqref{Eqn: Lying cost function}, \eqref{Eqn: network control cost functional}) while in protocol design it is measured through computational cost. 

Through these parallels, we see that introducing controls allows flexibility in constructing the problem as both the effect of the control and the cost function are chosen by the modeller. However, once these choices are made, one may hope that the resulting problem is well-posed, providing a unique optimal control whose form is not determined a priori. In contrast, protocol design problems may be more rigid in their setup but offer much greater freedom in their possible solutions, allowing modellers greater creativity in constructing protocols. A downside of this freedom is that it becomes more challenging to determine precisely when a `best' protocol has been found.

Perhaps reassuringly, we have seen throughout Section \ref{Section: Cts opinion dynamics} and Section \ref{Section: Binary Opinion Dynamics} that opinion manipulation is a challenging problem. Manipulations that are more direct (e.g. directly affecting individuals' opinions or interaction rules) are typically more effective in models than indirect manipulations (e.g. controlling network structure), even if the latter offers a greater number of levers to pull. However, such direct controls may be less effective in practice. Moreover, opinion manipulation problems become significantly more difficult when controls or protocols are limited to affect only a smaller number of individuals, edges, or interactions. Here, the typically non-linear and often stochastic nature of opinion formation models, which prevents straightforward predictions of the outcome of an intervention, can make such optimisation problems extremely challenging.

We would also like to acknowledge how manipulation and false claims, particularly spread on social media, disproportionally affect women. The recent rise in popularity of misogynistic figures such as Andrew Tate and Nick Fuentes is perpetuated by the use of manipulative and deceptive strategies. Much of the research that we have presented in this paper aims to answer the question of how best to control a population towards a desired opinion \cite{albi2014boltzmann, albi2014kinetic,glendinning2025what,nugent2024steering}. While we as researchers believe it is important to understand how deception is achieved in order to better protect marginalised and vulnerable groups, the optimal steering of a opinions is surely a question that many members of the so-called `manosphere' would be interested in. The ethics of modelling social systems are discussed in this recent paper \cite{aldana2025modeling}. Morally, we would like to encourage future authors in this field to perform research into how we might break these consensus control strategies and prevent influence of potentially dangerous manipulative figures. Work in this direction is ongoing, for example  D{\"u}ring et al. \cite{during2024breaking} present a method for breaking consensus by introducing a network to the problem, generalised to infinite dimensions as a graphon. 

\backmatter

\begin{appendices}

\section{Proof of Proposition \ref{Prop: Direct controllability}} \label{app: proof of controllability}

\begin{proof}
    Consider the control given by $u_i = x_d - x_i$. The dynamics then become
    \begin{align*}
        \dot{x}_i = \frac{1}{N}\sum_{j=1}^N\phi(x_j - x_i)(x_j - x_i) + x_d - x_i \,.
    \end{align*}
    We define the following functions 
    \begin{align*}
        \Phi(r) &= \int_0^r \phi(s) \, s \, ds \,, \\
        V(x) &= \frac{1}{2N} \sum_{i=1}^N \sum_{j=1}^N \Phi(x_j - x_i) + \frac{1}{2} \sum_{i=1}^N (x_d - x_i)^2 \,,
    \end{align*}
    then
    \begin{align*}
        \frac{dV}{dx_\ell} 
        &= \bigg( \frac{1}{2N} \sum_{j\neq\ell}^N -\phi(x_j - x_\ell)\,(x_j - x_\ell) + \frac{1}{2N} \sum_{i\neq\ell}^N \phi(x_\ell - x_i)\,(x_\ell - x_i) \bigg) - (x_d - x_\ell) \\
        &= - \bigg( \frac{1}{N} \sum_{j=1}^N \phi(x_j - x_\ell)\,(x_j - x_\ell)  + (x_d - x_\ell) \bigg)
    \end{align*}
    so we see that
    \begin{align*}
        \frac{dx}{dt} = -\nabla V (x) \,.
    \end{align*}
    Critical points of $V$ solve, for all $i=1,\dots,N$, 
    \begin{align*}
        \frac{1}{N}\sum_{j=1}^N\phi(x_j - x_i)(x_j - x_i) + x_d - x_i = 0 \,.
    \end{align*}
    Note that $x_i = x_d$ for all $i$ gives one solution to this equation. Let $x$ be any other solution, meaning there exists at least one $x_i \neq x_d$. Without loss of generality assume that the $x_i$'s are ordered so that $x_1\leq x_2 \leq \dots \leq x_N$. It must be the case that either $x_1 < x_d$ or $x_N > x_d$. 

    \textbf{Case 1:} $x_1 < x_d$
    
    As $x_1$ is the minimum opinion $x_j - x_1 \geq 0$ for all $j$. Moreover $\phi(x_j - x_1) \geq 0$ so 
    \begin{align*}
        0 = \frac{1}{N}\sum_{j=1}^N\phi(x_j - x_1)(x_j - x_1) + x_d - x_1 \geq x_d - x_1 > 0 \,,
    \end{align*}
    giving a contradiction. 

    \textbf{Case 2:} $x_N > x_d$
    
    As $x_N$ is the maximum opinion $x_j - x_N \leq 0$ for all $j$. Moreover $\phi(x_j - x_N) \geq 0$ so 
    \begin{align*}
        0 = \frac{1}{N}\sum_{j=1}^N\phi(x_j - x_N)(x_j - x_N) + x_d - x_N \leq x_d - x_N < 0 \,,
    \end{align*}
    again giving a contradiction. 

    In both cases we have a contradiction and so conclude that $x = x_d$ is the unique critical point of $V$. It is also clear that this a minimum of $V$. Thus by La Salle's invariance principle, $x_i \rightarrow x_d$ for all $i=1,\dots,N$. That is, we have convergence to consensus at $x_d$. 
\end{proof}

\section{Quasi-invariant limit of the Boltzmann description}\label{app:quasi-invariant-boltz}

We will present here the derivation of the Fokker-Planck equation from the Boltzmann equation in the case where agents' opinions are evolving according to the Hegelsmann-Krause model \eqref{eqn: HK-model}, following closely the derivation by Toscani in \cite{toscani2006kinetic}. The cases where we add control to our model, either by directly affecting opinions, incorporating leaders and followers or adding a liar follow similar derivations, explained in detail in \cite{albi2014kinetic, albi2014boltzmann} and \cite{glendinning2025what} respectively.

The first step to the derivation is to introduce the microscopic binary interactions that underpin the Boltzmann equation.

From a microscopic point of view, we describe a binary interaction between agents with opinions $(x, x_*)\in\mathcal{I}^2$ as the resulting change of opinion from said interaction, given by
\begin{align}\label{Eqn: appendix binary interactions}
\begin{split}
    x' &= x + \alpha \phi(x_* - x)(x_* - x) +  \theta D(x),\\
    x_*' &= x_* + \alpha \phi(x - x_*)(x - x_*) +\theta_{*}D(x_{*}),
\end{split}
\end{align}
where $(x', x_*')$ are the post-interaction opinions of $(x,x_*)$. Here, $\alpha\in(0,1/2)$ is a given constant, while $\theta, \theta_*$ are random variables with zero mean and variance $\sigma^2$. The function $D(\cdot)$ represents the diffusion of a given opinion.

We let $\mu(x,t)$ denote the distribution of opinion $x\in\mathcal{I}$ at time $t\geq0$. Standard methods of kinetic theory of binary interactions \cite{cercignani2013mathematical} are applied to recover a Boltzmann-type equation for time evolution of $\mu$ in weak form,
\begin{equation}\label{Eqn: appendix boltzmann}
    \frac{d}{dt}\int_{\mathcal{I}}\varphi(x)\mu(x,t)\, dx = (Q(\mu,\mu),\varphi),
\end{equation}
for $\varphi$ an arbitrary test function. Under the assumption that $D(x)$ and the noise $\theta,\theta'$ are suitable to preserve the bounds $x', x_*'\in\mathcal{I}$, we have that 
\begin{equation*}
    (Q(\mu,\mu),\varphi) = \left<\int_{\mathcal{I}^2}(\varphi(x') - \varphi(x))\mu(x,t)\mu(x_*,t)\, dx_*\, dx\right>,
\end{equation*}
where $\left<\cdot\right>$ denotes expectation with respect to $\theta,\theta_*$ and $x'$ is the post-interaction opinion of an agent with opinion $x$ after interaction with an agent of opinion $x_*$, given by \eqref{Eqn: appendix binary interactions}.

For the quasi-invariant limit, we rescale time $t$, the propensity strength $\alpha$ and the diffusion strength $\sigma^2$ simultaneously in order to maintain the memory of microscopic interactions at an asymptotic level. As in \cite{toscani2006kinetic}, we make the following scaling assumptions,
$$ \alpha = \varepsilon, \quad t = \frac{\tau}{\varepsilon}, \quad \sigma^2 = \varepsilon\tilde{\sigma}^2,$$
where $\varepsilon>0$. With this rescaling, equation (\ref{Eqn: appendix boltzmann}) becomes 
\begin{equation*}
    \frac{d}{d\tau}\int_{\mathcal{I}}\varphi(x)\mu(x,\tau)\, dx = \frac{1}{\varepsilon}\left<\int_{\mathcal{I}^2}(\varphi(x') - \varphi(x))\mu(x,\tau)\mu(x_*,\tau)\, dx_*\, dx\right>.
\end{equation*}
To recover the Fokker-Planck equation corresponding to the Boltzmann equation (\ref{Eqn: appendix boltzmann}), we are interested in the asymptotic behaviour of the equation as $\varepsilon\to 0$. To this end, we consider the second order Taylor expansion of $\varphi$ around $x$,
\begin{equation*}
    \varphi(x') -\varphi(x) = (x' - x)\varphi'(x) + \frac{1}{2}(x' - x)^2\varphi''(\tilde{x}),
\end{equation*}
where
$$ \tilde{x} = \vartheta x' + (1-\vartheta)x$$
for some $0\leq \vartheta\leq 1$. The rescaled binary dynamic (\ref{Eqn: appendix binary interactions}) can be written as
\begin{equation*}
    x' - x = \varepsilon\phi(x_* - x)(x_* - x) + \theta^{\varepsilon}D(x) + \mathcal{O}(\varepsilon^2),
\end{equation*}
where $\theta^{\varepsilon}$ is a random variable with mean 0 and variance $\varepsilon\tilde{\sigma}^2$. We hence have that the rescaled interaction integral $(Q(\mu, \mu), \varphi)$ can be written as
\begin{equation*}
    \int_{\mathcal{I}^2}\left[\phi(x_* - x)(x_* - x) \varphi'(x) + \frac{\tilde{\sigma}^2}{2}D(x)^2\varphi''(x)\right]\mu(x)\mu(x_*)\, dx_*\, dx + R(\varepsilon) + \mathcal{O}(\varepsilon),
\end{equation*}
where the remainder term $R(\varepsilon)$ is given by
\begin{equation*}
    R(\varepsilon) = \frac{1}{2\varepsilon}\left<\int_{\mathcal{I}^2}(x'-x)^2(\varphi''(\tilde{x})-\varphi''(x))\mu(x)\mu(x_*)\, dx_*\, dx\right>.
\end{equation*}
It can be shown rigorously, following closely from the derivation by Toscani \cite{toscani2006kinetic}, that $R(\varepsilon)\to 0$ as $\varepsilon\to 0$. Therefore, we have the following limiting equation for (\ref{Eqn: appendix boltzmann}),
\begin{equation*}
    \frac{d}{d\tau}\int_{\mathcal{I}}\varphi(x)\mu(x,\tau)\, dx = \int_{\mathcal{I}^2}\phi(x_* - x)(x_* - x)\varphi'(x)\, \mu(x)\mu(x_*)\, dx_* \, dx+ \frac{\tilde{\sigma}^2}{2}\int_{\mathcal{I}}D(x)^2\varphi''(x) \mu(x)\, dx.
\end{equation*}
Integrating the right-hand-side by parts and applying the variational lemma, since $\varphi$ is arbitrary, we have the resulting Fokker-Planck equation \eqref{eqn: mean-field-pde} given in Section \ref{Section: Introduction continuous models},
\begin{equation*}
    \partial_{\tau}\mu(x,\tau) + \partial_x\left(\int_{\mathcal{I}}\phi(x_*-x)(x_*-x)\mu(x_*,\tau)dx_*\mu(x,\tau)\right) = \frac{\tilde{\sigma}^2}{2}\partial_{xx}(D(x)^2\mu(x,\tau)).
\end{equation*}

\end{appendices}

\bibliography{bibliography}

\end{document}